# *Triangoli, Icosaedri e Cupole Geodetiche*

Giuseppe Conti* Raffaella Paoletti**

\* Giuseppe Conti, Dipartimento di Matematica DIMAI, Università di Firenze. Italia. giuseppe.conti@unifi.it
\*\* Raffaella Paoletti, membro del GNSAGA, Dipartimento di Matematica DIMAI, Università di Firenze, Italia. raffaella.paoletti@unifi.it

***Sommario*** *Le cupole geodetiche, poliedri convessi di forma quasi sferica o parti di essi, sono state oggetto di grande attenzione nel ventennio compreso fra le metà degli anni Cinquanta e Settanta del secolo scorso, soprattutto grazie a Richard Buckminster Fuller. Dopo un boom edilizio, specialmente negli Stati Uniti, il loro interesse costruttivo è declinato ma le loro caratteristiche geometriche, studiate da vari matematici, hanno inaspettatamente trovato applicazioni in altri campi delle Scienze (Biologia, Chimica) fin dalla metà degli anni Cinquanta e sono tuttora alla base di modelli, oggetto di ricerche attuali. Anche in campo ingegneristico-architettonico, col rinato interesse per le "grandi" strutture reticolari, esse vengono analizzate e a volte prese come spunto. In questo articolo, dopo aver riassunto la storia della loro ideazione e delle varie applicazioni nelle Scienze, ne descriviamo le principali caratteristiche geometriche ed i procedimenti geometrici più usati per la loro costruzione.*

***Abstract*** *Geodesic domes, convex polyhedrons with almost spherical shape or parts of them, were the subject of great attention in the twenty years between the mid-1950s and the 1970s, especially thanks to Richard Buckminster Fuller. After a building boom, mostly in the United States, their construction interest declined but their geometric characteristics, studied by various mathematicians, have unexpectedly found applications in other fields of science (Biology, Chemistry) since the mid-1950s and are still underlying models, the subject of current research. Even in the engineering-architectural field, with the revived interest in "large" reticular structures, they are analyzed and sometimes taken as inspiration. In this article, after summarizing the history of their conception and the various applications in science, we describe their main geometric characteristics and the geometric procedures most used for their construction.*

## *1 - Motivazioni*

Il triangolo, poligono col numero minimo di vertici e lati, non ha diagonali e questo lo rende una figura **rigida** o **indeformabile**: tre aste unite a triangolo formano una struttura che non si deforma anche se le connessioni alle estremità non sono completamente rigide. Una poligonale con un numero maggiore di lati costruita in modo analogo è deformabile, ma diventa rigida inserendo ulteriori aste per "triangolarla". Questa proprietà è stata sfruttata fin dall'antichità nell'ambito della statica edilizia, in particolare per le coperture; a partire dal quarto secolo d.C., le travi orizzontali delle coperture vengono sostituite da semplici capriate in legno formate da due *puntoni*; successivamente, per aumentarne la resistenza, viene introdotta la *catena* (elemento orizzontale) per passare poi a capriate più articolate con *monaco* ed eventuali *saettoni* per una triangolazione più fitta. Nel Medioevo queste coperture vengono usate per spazi ampi come quelli dei monasteri e delle cattedrali, in alcuni casi con risultati considerati ancora oggi notevoli.

Intorno alla metà dell'Ottocento si è fatto ricorso a strutture reticolari e, successivamente, alle *travature reticolari*: una sorta di "sandwich" formato da due sistemi continui paralleli (*correnti*) che racchiudono una struttura reticolare, nella maggior parte dei casi con reticolo triangolare. Queste strutture hanno permesso di coprire grandi spazi e sono (state) utilizzate, ad esempio, per coperture industriali oppure ponti e passerelle; il loro ampio uso è dovuto alla loro resistenza a forti carichi pur avendo un peso ed un costo ridotto rispetto a molte altre soluzioni. Fra i numerosi esempi di strutture reticolari ricordiamo gli iconici *Forth Rail Bridge* di Edimburgo (lungo 2,5km, opera di Fowler e Baker, 1890) e la *Torre Eiffel* (Parigi, 1889).

Intorno alla metà del Novecento vengono sviluppati altri modelli di "triangolazione" ed uno dei più importanti è sicuramente quello delle **cupole geodetiche**: (porzioni di) poliedri convessi inscritti in una sfera i cui spigoli (rettilinei) giacciono, approssimativamente, lungo circonferenze massime (le *linee geodetiche* della sfera). Le facce di queste cupole sono, in genere, dei triangoli. La struttura ottenuta tende a ripartire gli sforzi in modo uniforme sulle aste, garantendo rigidezza e staticità a fronte di un peso contenuto.

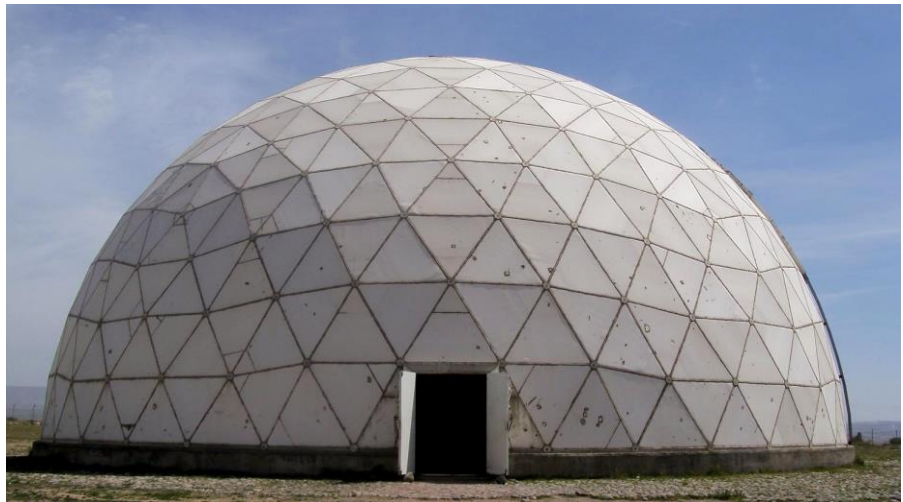

*Figura 1. Osservatorio di Maragheh (Azerbaijan).*
*Foto di: Elmju - Opera propria, CC BY-SA 3.0*
*https://commons.wikimedia.org/w/index.php?curid=10027825*

Dato che la sfera permette di racchiudere il massimo volume a parità di superficie, la copertura a forma (di porzione) sferica permette di usare minor materiale a parità di volume racchiuso. Ad esempio, una casa cubica con lato di misura $l$ racchiude un volume $l^3$ ed ha una superficie esterna $5l^2$; una semisfera di raggio $\sqrt{\frac{5}{2\pi}}\, l$ ha la stessa superficie esterna ma racchiude un volume $\frac{5\sqrt{5}}{3\sqrt{2\pi}}\, l^3$, circa una volta e mezzo quello del cubo[1].

Progettata dallo studio *Synergetics* di Richard Buckminster Fuller, nel 1958 viene costruita a Baton Rouge (Louisiana) per la Union Tank Car Company una cupola di 118mt di diametro; all'epoca è la copertura più ampia di uno spazio completamente libero all'interno[2].

Prima di vedere nel dettaglio la genesi e la costruzione di queste cupole, osserviamo alcune delle "grandi" costruzioni contemporanee realizzate con strutture reticolari triangolari.

*Tacoma Dome*, cupola non geodetica; inaugurata nel 1983, con i suoi 160mt di diametro è la cupola pubblica più grande al mondo. Nel 1983 era anche la copertura in legno più grande al mondo.

*Queen Elizabeth II Great Court*, copertura della corte del British Museum inaugurata nel 2000 e progettata dallo studio *Foster & partners*; formata da 3212 pannelli di vetro triangolari con forme diverse, è la piazza coperta più grande d'Europa.

Sia il grattacielo *The Gherkin* (Londra, 2004) che la *Hearst Tower* (NewYork, 2006), progettati da Norman Foster, si caratterizzano per una struttura a maglie triangolari detta *diagrid* (contrazione dei termini *diagonal* e *grid*) che si eleva per oltre 180mt. Queste costruzioni sono i primi esempi di applicazione della *diagrid* ad edifici molto alti.

Nel 2016 la *Carlo Ratti Associati* presenta il progetto "avveniristico" *The Mile*: si tratta di un giardino verticale alto 1609mt (un miglio appunto) basato su una struttura reticolare triangolata, larga 20mt, che si avvolge a forma di elicoide.

Ovviamente la statica delle opere menzionate è estremamente complessa, ma quelle a cupola geodetica devono sostanzialmente la loro stabilità ad un teorema di matematica: il *Teorema di rigidità* di *Augustin Louis Cauchy*, il quale afferma che un poliedro convesso è indeformabile (come struttura nello spazio) se le sue facce (che sono figure piane) sono rigide[3].

---

[1] Per completezza osserviamo anche che, ovviamente, la superficie di base della semisfera è 5/2 più grande di quella del cubo. Negli Stati Uniti di metà Novecento, dove vengono costruite la maggior parte delle cupole, questo parametro è secondario rispetto al risparmio di materiale.

[2] Nel 1990 l'edificio viene venduto alla Kansas City Southern Railway che poi lo demolirà a partire dal 2007, nonostante le numerose proteste popolari. Nel 1961, sempre per la Union Tank Car Company, viene costruita una cupola uguale a questa a Wood River, Illinois, ancora esistente sotto altra proprietà.

Intorno al 1940, **Aleksandr Danilovič Aleksandrov** dimostrò che tutte le superfici poliedriche convesse triangolate sono rigide se non contengono vertici all'interno delle facce naturali piatte[4].

La questione se tutti i poliedri, anche quelli non convessi, fossero rigidi o meno è rimasta aperta a lungo[5]. Nel 1975 **Herman Gluck** ha dimostrato che quasi tutti i poliedri sono rigidi[6].

A quel tempo, molti credevano che tutti i poliedri fossero rigidi (seguendo una congettura di Eulero del 1766), ma nel 1977 **Robert Connelly** sorprese la comunità scientifica costruendo un poliedro non rigido nello spazio tridimensionale; tale poliedro è a facce triangolari ma, naturalmente, non è convesso[7]. Notiamo che, per il risultato di Gluck, si trattava di un esempio piuttosto raro tra i poliedri.

## 2 - Richiami sui poliedri

Il termine **poliedro** deriva dalle parole greche *poly* (molti) e *hedra* (sede) ma, nonostante questi solidi siano stati studiati fin dall'antichità (anche civiltà pre-elleniste come quelle egizia, babilonese e cinese, oltre agli stessi greci, hanno calcolato il volume di un tronco di piramide), una definizione matematica di poliedro è stata data solo in secoli più recenti. Senza entrare nel merito delle "varie" definizioni di poliedro che sono attualmente usate, definiamo **poliedro convesso** la regione di spazio convessa delimitata da un numero finito di facce poligonali in modo che valgano le seguenti tre condizioni:
1- due facce sono disgiunte oppure si intersecano in uno spigolo o in un vertice,
2- ogni spigolo appartiene esattamente a due facce,
3- due facce adiacenti non sono complanari.

La porzione di spazio delimitata dai due semipiani che contengono due facce adiacenti viene chiamata *diedro*, mentre tutte le facce che si intersecano in uno stesso spigolo determinano un *angolo solido* o *angoloide*.

Si chiama **grado (o valenza) di un vertice** il numero degli spigoli del poliedro che hanno quel vertice come estremo. Osserviamo che lo stesso numero rappresenta anche il numero delle facce a cui appartiene quel vertice.

Un **poliedro** convesso si dice **regolare** se tutte le sue facce sono poligoni regolari congruenti tra loro e tutti i vertici hanno lo stesso grado.

Come è noto, si possono costruire poligoni regolari con un numero *n* di lati per ogni intero $n \geq 3$, ma esistono solo 5 poliedri regolari, detti **solidi platonici** perché Platone ne parla esplicitamente nella sua opera *Il Timeo* [Fr].

Tetraedro, cubo e dodecaedro erano già noti alla scuola pitagorica ed in seguito i greci introdussero l'ottaedro e l'icosaedro [Br]; nel 410 a.C. Teeteto dimostrò che i poliedri regolari sono soltanto questi cinque [Br], mentre Platone legò questi solidi agli elementi della Natura: il cubo alla Terra (stabilità), il tetraedro al Fuoco (ha il minor numero di lati), l'ottaedro all'Aria (il più mobile dei solidi), l'icosaedro all'Acqua (ha il maggior numero di lati) e il dodecaedro all'Universo (ha 12 facce come le case dello zodiaco). Si ipotizza che gli studi degli antichi greci sui poliedri siano nati dall'osservazione della natura: infatti, cubo, ottaedro e dodecaedro (quest'ultimo non regolare) sono tra le forme dei cristalli di fluorite e pirite, mentre altri solidi (fra cui l'icosaedro) sono visibili negli scheletri di alcuni animali marini come i Radiolari [DA][Cr].

Non è chiaro se i solidi platonici siano stati oggetto di studio di popolazioni precedenti: ad esempio, sono stati ritrovati in Scozia (ma anche in Inghilterra e Irlanda) vari manufatti in pietra risalenti al neolitico (circa 3000-2500 a.C.) che hanno la forma di sfere intagliate in modo da far risultare vari "pomelli" sporgenti.

---

[3] Vedi: Cauchy A. (1813), *Sur les polygones et polyèdres*, J. Éc. Polytech., IX, pp. 87-99.
[4] Vedi: Aleksandrov, A. D. (Traduzione di Zalgaller 2006), *Convex Polyhedra*, Springer Monographs in Mathematics, Springer.
[5] Il problema della rigidità (globale e infinitesimale) dei poliedri e delle strutture articolate è stato oggetto di studio di diverse branche della Matematica.
[6] Vedi: Gluck H. (1975), *Almost all simply connected closed surfaces are rigid,* Lecture Notes in Math. 438, Geometric Topology, Springer-Verlag, pp. 225-239.
[7] Vedi: Connelly R. (1977), *A counterexample to the rigidity conjecture for polyhedra,* Publications mathématiques de l'I.H.É.S., Tome 47, pp. 333-338.

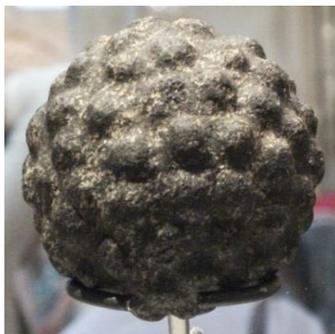

*Figura 2. Foto di: Johnbod [licenza CC]*
*https://en.wikipedia.org/wiki/Carved_stone_balls#/media/File:Room_51,_British_MuseumDSCF6620.jpg*

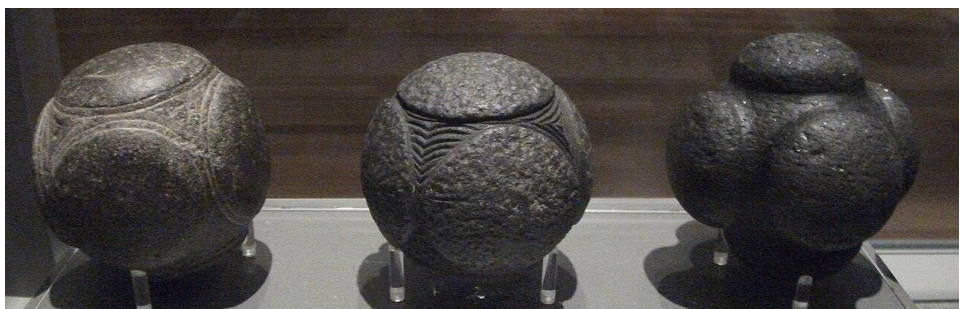

*Figura 3. Autore: Johnbod [licenza CC]*
*https://en.wikipedia.org/wiki/Carved_stone_balls#/media/File:Kelvingrove_Art_Gallery_and_MuseumDSCF0239_11.JPG*

Il loro uso è ancora incerto e dibattuto fra gli studiosi; il numero di "pomelli" è variabile da 3 a 160 ed alcuni sono decorati con incisioni. Potrebbero essere stati ideati per motivi mistico-religiosi oppure decorativi o anche dall'osservazione della natura. Benché alcuni abbiano una forma molto simile ai solidi platonici, gli studiosi ritengono che non siano frutto di studi matematici ma siano stati realizzati così allo scopo di distribuire in maniera regolare i vari pomelli.

Sono invece dichiaratamente a forma di dodecaedro e icosaedro numerosi manufatti bronzei di epoca romana, datati dal II al IV secolo d.C. e rinvenuti soprattutto in Francia e Germania; si suppone siano stati considerati oggetti di valore ma il loro impiego non è chiaro, anche perché i manufatti hanno dimensioni diverse (da 4 a 10cm). Interessante una recente ipotesi [Sp] che vede il dodecaedro usato come telemetro (per gli appezzamenti agricoli).

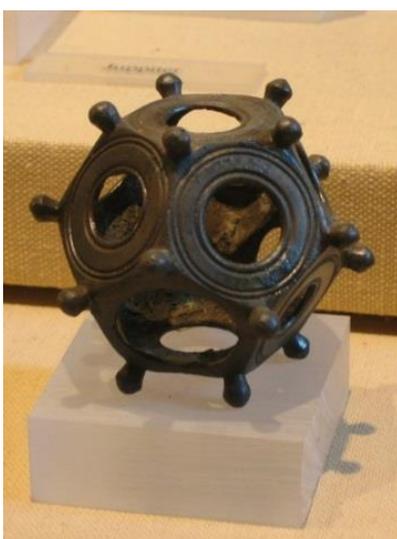

*Figura 4. Foto di: User Itub [licenza CC]*
*https://it.m.wikipedia.org/wiki/File:Roman_dodecahedron.jpg*

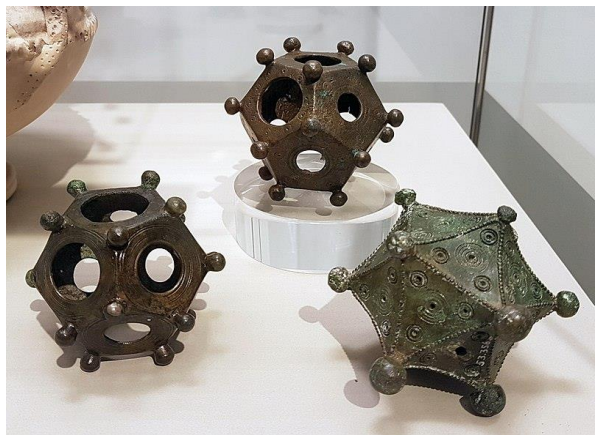
*Figura 5. Foto di: User Kleon3   [licenza CC]*
*https://commons.wikimedia.org/wiki/File:2018_Rheinisches_Landesmuseum_Bonn,_Dodekaeder_%26_Ikosaeder.jpg*

Nel IV secolo d.C. Pappo d'Alessandria elencò e studiò una nuova famiglia di poliedri, che attribuì ad Archimede [Br] e che sono chiamati oggi *poliedri archimedei*. Questi sono dei poliedri convessi in cui le facce sono dei poligoni regolari (dunque gli spigoli del solido sono tutti congruenti), tutti i vertici hanno lo stesso grado, ma in modo che il solido non sia né un poliedro platonico né un prisma o un antiprisma. Per ottenere un poliedro archimedeo si devono, quindi, usare almeno due tipi diversi di poligoni regolari; questi solidi sono quelli più regolari dopo quelli platonici e per questo sono detti anche *poliedri semiregolari*. Ricordiamo che *prismi* e *antiprismi* sono poliedri formati da due poligoni con $n$ lati (le basi), uguali e paralleli, collegati fra loro da $n$ parallelogrammi (prismi) o da triangoli isosceli (antiprismi; in questo caso i triangoli sono $2n$ e ogni triangolo ha due vertici su una delle basi ed il terzo vertice sull'altra base). Se i poligoni delle basi sono regolari e i parallelogrammi del prisma sono quadrati, oppure i triangoli dell'antiprisma sono equilateri, si hanno due famiglie infinite di solidi "quasi regolari" che non hanno però molte simmetrie.

I solidi archimedei sono 13 e sono tutti inscrivibili in una sfera, ma non circoscrivibili [C-R]. Cinque di questi possono essere ottenuti per *troncamento* (ovvero tagliando i vertici) dei solidi platonici in modo che gli spigoli del solido così ottenuto siano tutti uguali alla terza parte nei poliedri a facce triangolari, poco meno negli altri due, degli spigoli del poliedro di partenza (vedi [Ba]). Ne è un esempio l'*icosaedro troncato*, con facce pentagonali ed esagonali, utilizzato dalla ditta Adidas per modellare il famoso pallone da calcio *Telstar* alla fine degli anni 1960 [C-T-C].

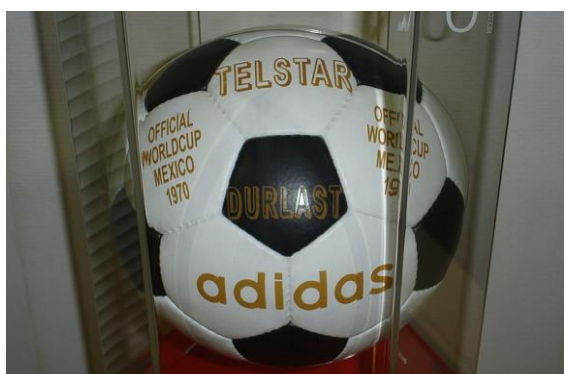
*Figura 6. Foto di User Zac Allan*
*https://commons.wikimedia.org/wiki/File:Adidas_Telstar_Mexico_1970_Official_ball.jpg*

Altri due poliedri archimedei si ottengono eseguendo il troncamento a metà degli spigoli dei solidi platonici, in modo che gli spigoli del solido ottenuto siano tutti uguali alla metà degli spigoli del poliedro di partenza. Questa costruzione dà luogo solo a due poliedri archimedei perché sia dal cubo che dall'ottaedro si ottiene il cubottaedro, sia dal dodecaedro che dall'icosaedro si ottiene l'icosidodecaedro, mentre dal tetraedro si ottiene l'ottaedro che è un solido platonico[8]. Per ottenere gli altri sei poliedri archimedei si hanno differenti costruzioni. Notiamo che fra i poliedri archimedei soltanto il *cubo camuso* e il *dodecaedro camuso* sono *chirali*[9].

---

[8] Vedremo fra poco che cubo e ottaedro, dodecaedro e icosaedro sono due coppie di solidi duali.

Un poliedro convesso, le cui facce sono tutti poligoni regolari che non è né un solido platonico o archimedeo, né un prisma o antiprisma, è detto **solido di Johnson**. L'esempio più semplice di questo tipo è la piramide a base quadrata avente come facce laterali 4 triangoli equilateri. Furono scoperti da Norman Johnson [Jo] che ne elencò 92; nel 1967 Victor Zalgaller dimostrò che effettivamente sono esattamente 92.

Introduciamo adesso il concetto di **dualità** fra poliedri, di cui esistono varie definizioni. Seguendo la costruzione illustrata in [C-R], basata sulla *polarità* indotta da una sfera, possiamo associare ad ogni punto dello spazio un piano (detto **piano polare**); in tal caso il punto è detto **polo** del piano. Se il punto appartiene alla sfera, allora il suo piano polare è il piano tangente alla sfera in quel punto[10]. Si ottiene così una dualità fra punti e piani che soddisfa la legge di reciprocità: se un punto Q appartiene al piano polare di un punto P, allora il piano polare di Q contiene P. In questa corrispondenza, ad una retta passante per due punti P e Q viene associata la retta intersezione dei piani polari di P e Q. La polarità trasforma un poliedro *P* in un altro poliedro *P\**, detto **duale** di *P*, in modo che ad ogni vertice/spigolo/faccia di *P* corrisponda faccia/spigolo/vertice di *P\** con le seguenti inclusioni:

1- se una faccia F di *P* contiene uno spigolo S, allora il vertice di *P\** trasformato di F è contenuto nello spigolo trasformato di S;
2- se V è un vertice di *P* contenuto in uno spigolo S, allora la faccia di *P\** trasformata di V contiene lo spigolo trasformato di S.

Il duale del duale di un poliedro è il poliedro iniziale.

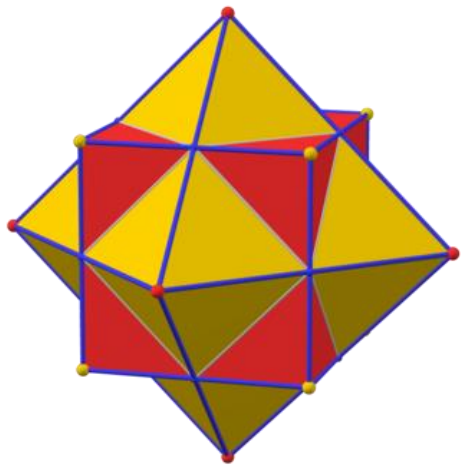

*Figura 7. Autore T.Piesk  Il duale di un cubo (in rosso) è l'ottaedro (in giallo). La costruzione del duale in questo caso è stata fatta rispetto alla sfera tangente a tutti gli spigoli del cubo (intersfera).*
*https://commons.wikimedia.org/wiki/File:Polyhedron_pair_6-8.png*

Poiché siamo interessati a poliedri convessi inscrivibili in una sfera, possiamo definire la polarità rispetto alla sfera circoscritta al poliedro, ottenendo in questo modo un metodo più semplice per costruire il poliedro duale. Infatti, in questo caso il poliedro duale è dato dall'intersezione dei piani tangenti alla sfera nei punti che sono vertici del poliedro di partenza; viceversa, il duale di un poliedro circoscrittibile alla sfera ha come vertici i punti di tangenza con la sfera inscritta.

Si può ottenere il duale dei poliedri platonici anche con una costruzione più semplice: consideriamo i centri delle facce di questi solidi e congiungiamo con un segmento due centri che appartengono a facce adiacenti. Si ottiene così un poliedro simile a quello ottenuto con la costruzione mediante la sfera circoscritta, quindi tale poliedro è il duale del solido di partenza.

---

[9] Un solido si dice **chirale** (dalla parola greca χείρ che significa "mano") se non è sovrapponibile, tramite rotazioni e traslazioni, alla sua immagine speculare. Esempi di oggetti chirali sono le nostre mani e i nostri piedi. La definizione di figura chirale può essere data anche nel piano; ad esempio, un triangolo scaleno è una figura piana chirale.
[10] La polarità può essere definita anche attraverso l'inversione per raggi reciproci determinata dalla sfera: se P' è il corrispondente nell'inversione di un punto P, allora il piano polare è il piano per P' avente come vettore normale PP'.

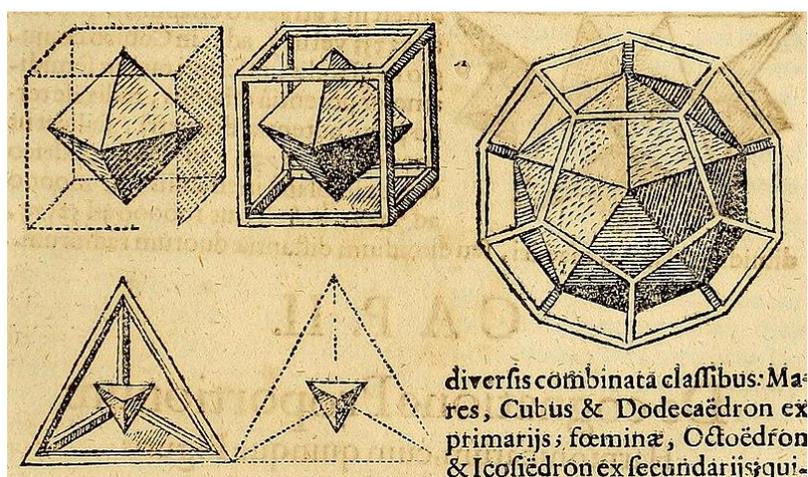
*Figura 8. Tratta da Harmonices Mundi di J.Kepler*
*https://archive.org/details/ioanniskepplerih00kepl/page/180/mode/2up?view=theater*

Come si vede dalla Figura 8, il duale di un solido platonico è ancora un solido platonico: il duale del cubo è l'ottaedro, quello dell'icosaedro è il dodecaedro e il tetraedro è autoduale.

I duali dei poliedri archimedei sono i ***solidi*** (detti anche ***poliedri***) ***di Catalan***, descritti nel 1865 dal matematico belga Eugène Charles Catalan. Poiché i poliedri archimedei possono essere inscritti in una sfera, i loro duali (i poliedri di Catalan) hanno le facce tangenti alla sfera; viceversa, i duali dei poliedri di Catalan (i poliedri archimedei) hanno per vertici i punti di tangenza delle facce alla sfera inscritta al poliedro di Catalan.

Un'altra importante classe di poliedri è data dai ***poliedri di Goldberg***, che prendono il nome dall'ingegnere (appassionato di matematica) Michael Goldberg, che li descrisse nel 1937. Questi sono poliedri convessi con le seguenti caratteristiche: tutti i vertici hanno grado tre, le facce sono solo pentagoni o esagoni, possiedono tutte le simmetrie rotazionali dell'icosaedro.

Studiando problemi isoperimetrici per poliedri, Goldberg introdusse una classe di poliedri con molte simmetrie, che chiamò *medial polyhedra* [Go1, Go2]; quello con 12 facce era il dodecaedro regolare, quelli con un numero di facce maggiore erano ottenuti componendo 12 pentagoni (regolari) con esagoni. Questi solidi hanno tutte le 60 simmetrie rotazionali dell'icosaedro, che sono: l'applicazione identica, 12 rotazioni di 72°, 12 rotazioni di 144°, 20 rotazioni di 120°, 15 rotazioni di 180°; a causa di queste simmetrie, i pentagoni che appaiono nei poliedri di Goldberg sono tutti regolari.

In letteratura si trovano varie definizioni dei poliedri di Goldberg, che differiscono fra loro per la regolarità richiesta alle facce esagonali; in questo articolo, considereremo solo i poliedri di Goldberg aventi esagoni equilateri (ma non necessariamente regolari); quindi, i poliedri di Goldberg che tratteremo hanno tutti gli spigoli congruenti fra loro.

Il numero fisso 12 dei pentagoni in questi solidi, come già calcolato da Goldberg, è una conseguenza della ***formula di Eulero: $V - S + F = 2$***. Infatti, indichiamo con $F_5$ e $F_6$ il numero delle facce pentagonali ed esagonali rispettivamente; poiché ogni spigolo è comune a due facce e in ogni vertice concorrono tre facce, si ha:

$$S = (5F_5 + 6F_6)/2, \quad V = (5F_5 + 6F_6)/3$$

Sostituendo tali valori nella formula di Eulero, si ottiene:

$$10F_5 + 12F_6 - 15F_5 - 18F_6 + 6F_5 + 6F_6 = 12,$$

da cui $F_5 = 12$.

Osserviamo che dalla precedente relazione segue anche: $V = 20 + 2F_6$. Quindi, un poliedro di Goldberg ha almeno 20 vertici e, se ne ha esattamente 20, allora non ha facce esagonali e, di conseguenza, è un dodecaedro.

Se $V = 60$ si ottiene invece l'icosaedro troncato; in tal caso $F_6 = 20$.

## 3 - Cupole geodetiche

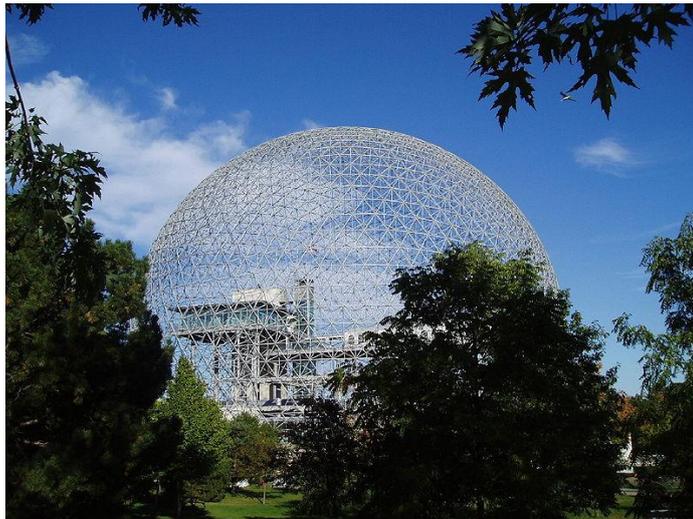

*Figura 9. Biosfera Montreal   Foto di Eberhard von Nellenburg [licenza CC]*
*https://commons.wikimedia.org/wiki/File:Mtl._Biosphere_in_Sept._2004.jpg*

In letteratura, col termine *cupola geodetica* si indica di solito una griglia (semi)sferica a maglie triangolari le cui aste giacciono, approssimativamente, sulle geodetiche (ovvero i cerchi massimi) della sfera. In questo articolo, chiameremo *sfera geodetica* (o *poliedro geodetico*) un poliedro topologicamente equivalente ad una sfera, a facce triangolari, inscritto in una sfera, i cui vertici hanno grado 5 o 6 e i cui lati giacciono approssimativamente sulle geodetiche della sfera. Chiameremo *cupola geodetica* una parte di questo poliedro, ad esempio quella corrispondente ad una semisfera o ad una generica calotta sferica.

Di solito le cupole geodetiche realizzate in architettura hanno facce triangolari proprio per il citato Teorema di Cauchy: in questo modo la struttura risultante è stabile anche se le aste dell'intelaiatura non hanno connessioni completamente rigide agli estremi poiché, comunque, formano dei poligoni rigidi. Per costruire cupole geodetiche con facce non triangolari occorre rendere rigide tali facce.

In generale le aste sono realizzate con profili tubolari di sezione circolare e sono dotate alle loro estremità di appositi innesti con dei fori che permettono tramite un bullone di unirle ai vari nodi, componendo in tal modo la struttura geometrica. Le aste sono incernierate internamente sia nel piano meridiano che in quello tangente in modo tale che l'unica sollecitazione che si trasferisce da un'asta all'altra sia quella normale. Notiamo che i nodi delle cupole geodetiche possono essere schematizzati con cerniere oppure con incastri.

Nelle cupole geodetiche, gli spigoli distribuiscono gli sforzi locali sull'intera struttura; queste cupole sono concepite in modo che eventuali azioni esterne (forze), ad esempio quelle dovute ad un telo di copertura, vengano applicate solo nei nodi. Inoltre, a parità di forma, dimensioni complessive, materiale utilizzato, la loro resistenza aumenta se aumenta il numero di triangoli che le compongono; quindi, strutture più grandi sono anche più resistenti. Le deformazioni di queste strutture sono trascurabili poiché sono sollecitate principalmente a trazione oppure a compressione; tuttavia, occorre osservare che il comportamento delle forze di trazione e di compressione nelle differenti strutture geodetiche non è stato ancora completamente compreso.

Come vedremo in dettaglio più avanti, uno dei metodi più usati per costruire sfere geodetiche a facce triangolari consiste nel suddividere ogni faccia di un **icosaedro** regolare (il poliedro regolare che meglio approssima una sfera) in triangoli equilateri più piccoli e tutti uguali fra loro; successivamente dal centro della sfera circoscritta all'icosaedro si proiettano tutti i vertici così ottenuti sulla sfera stessa. In ogni vertice del poliedro ottenuto concorrono 5 o 6 spigoli; perciò, i vertici sono di grado 5 o 6.

I poliedri duali di queste sfere sono poliedri a facce pentagonali o esagonali, hanno tre facce concorrenti in ogni vertice e, per come sono stati costruiti, hanno le simmetrie rotazionali di un icosaedro, quindi sono poliedri di Goldberg.

Nelle cupole geodetiche la proiezione ortogonale del centro della sfera su ciascuna faccia (triangolare) del poliedro è il circocentro del triangolo. Questa proprietà è molto utile nella costruzione della cupola e si dimostra applicando il *Teorema delle tre perpendicolari.* Infatti, siano A, B, C i vertici di una faccia *F*, O il centro della sfera circoscritta alla cupola e K il piede della perpendicolare condotta da O al piano in cui giace *F*. Nel piano contenente *F*, conduciamo da K la perpendicolare al lato AB e indichiamo con H il suo punto di intersezione con AB. Il piano contenente i punti O, K e H è quindi, per il teorema delle tre perpendicolari, ortogonale ad AB; in particolare, OH è ortogonale ad AB. Poiché il triangolo AOB è isoscele rispetto alla base AB e OH è l'altezza relativa a tale base, si ha che H è il punto medio di AB; quindi, KH è l'asse di AB nel piano di *F*. Analoghe considerazioni valgono per i lati AC e BC, quindi K è il circocentro di ABC[11].

Un altro esempio interessante di sfera geodetica si ottiene a partire dal **dodecaedro** regolare con una costruzione differente: su ogni faccia si costruisce una piramide (a base pentagonale) il cui vertice è dato dall'intersezione della sfera circoscritta al dodecaedro con la retta passante per il centro della sfera e perpendicolare alla faccia pentagonale; con questa costruzione si ha che tale piramide è retta. Il poliedro che si ottiene in questo modo ha come facce tutti triangoli isosceli congruenti e fa parte della famiglia dei **pentacisdodecaedri**[12]. Notiamo che i vertici delle piramidi costruite sono i vertici dell'icosaedro duale del dodecaedro di partenza.

Un poliedro di questo tipo fu progettato da *Michelangelo Buonarroti* per il coronamento della lanterna della cupola della Sagrestia Nuova della chiesa di San Lorenzo a Firenze. Quasi sicuramente questa *palla*, come la chiama Michelangelo, è inscrivibile in una sfera. L'esecuzione di questa sfera geodetica, formata da due pezzi non saldati fra loro, fu opera dell'orafo e scultore fiorentino *Giovanni di Baldassarre* (1460 - 1536), specializzato nelle opere di metallo e soprannominato il **Piloto.**

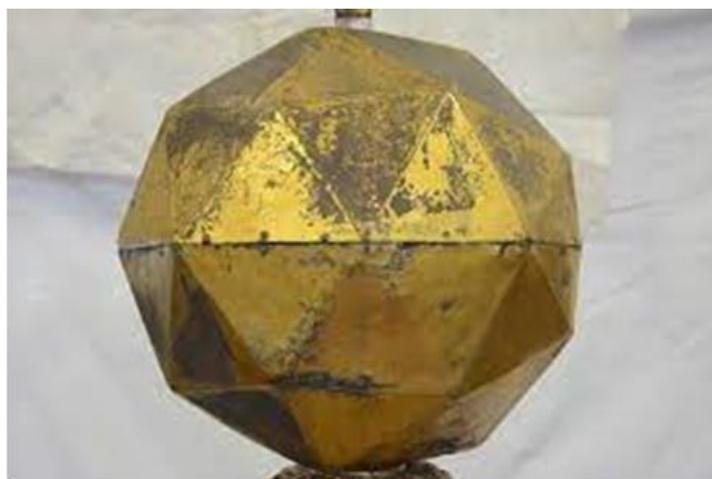

*Figura 10. Il poliedro posto a coronamento della lanterna della Sagrestia Nuova di San Lorenzo a Firenze*
**http://met.provincia.fi.it/public/images/20130117142819681.jpg**

Una curiosità: le facce di un pentacisdodecaedro sono ovviamente 60 ma, per motivi che ci sfuggono, Giorgio Vasari afferma che sono 72 [Va, p. 1223]!

Un pentacisdodecaedro è stato disegnato da Leonardo da Vinci nel libro di Luca Pacioli *Divina Proportione*; tuttavia, questo poliedro è un po' diverso da quello di Michelangelo perché nel poliedro di Leonardo le facce delle piramidi sono triangoli equilateri; quindi, il poliedro ottenuto non è convesso e non è inscrivibile in una sfera.

---

[11] Notiamo che in realtà questa proprietà vale per qualunque poliedro a facce triangolari (acutangole) inscritto in una sfera.
[12] I pentacisdodecaedri sono poliedri ottenuti costruendo su ogni faccia di un dodecaedro regolare una piramide retta avente per base la suddetta faccia. Le facce di tutte queste piramidi sono triangoli isosceli congruenti fra loro.

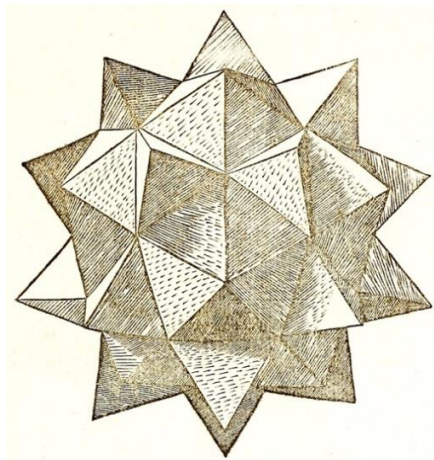

*Figura 11. Il pentacisdodecaedro di Leonardo da Vinci*
*https://commons.wikimedia.org/wiki/File:De_divina_proportione_-_Illustration_08.jpg*

Esiste anche il pentacisdodecaedro circoscrivibile ad una sfera: si tratta del poliedro di Catalan duale del solido archimedeo icosaedro troncato. Le facce del pentacisdodecaedro di Catalan sono triangoli isosceli tutti uguali fra loro; si dimostra che in ogni triangolo il rapporto fra il lato obliquo e la base è circa 0,887057998, mentre l'angolo al vertice è circa 68° 37' 07".

Anche nel pentacisdodecaedro inscritto in una sfera le facce sono triangoli isosceli tutti uguali fra loro, ma in questo poliedro il rapporto fra il lato obliquo e la base di ogni triangolo è circa 0,897999085 e l'angolo al vertice è circa 67°40' 07". Poiché le differenze fra i due pentacisdodecaedri (quello inscritto e quello circoscritto) sono minime, resta difficile stabilire con certezza quale sia quello di Michelangelo. Noi pensiamo che, per questioni realizzative, l'artista abbia pensato a quello inscritto.

Con un procedimento analogo a quello appena visto per costruire il pentacisdodecaedro inscritto nella sfera, si può costruire una sfera geodetica a partire da un icosaedro troncato: si proiettano i centri delle facce (che sono pentagoni o esagoni regolari) dal centro della sfera circoscritta al solido sulla sfera stessa e si costruiscono delle piramidi aventi per vertice i punti proiettati e per base la faccia di partenza. In questo caso si ottengono 12 piramidi a base pentagonale e 20 a base esagonale; quindi, la sfera geodetica ottenuta ha 270 spigoli (di 3 lunghezze diverse), 180 facce (triangoli isosceli di due tipi) e 92 vertici (di due tipi)[13].

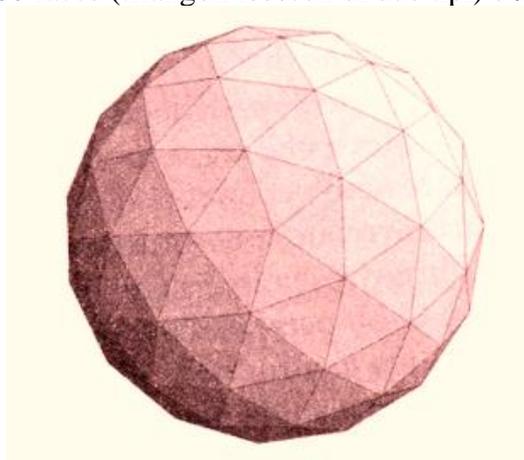

*Figura 12. Sfera geodetica ottenuta dall'innalzamento dell'icosaedro troncato*

Il procedimento costruttivo appena visto, detto anche **twinning (gemmazione)**, può essere esteso ad un qualunque poliedro inscrivibile in una sfera ed avente come facce dei poligoni regolari, non triangolari. Tuttavia, partendo da solidi diversi da quelli considerati adesso, le lunghezze dei lati del poliedro e degli spigoli delle piramidi possono essere molto diverse e dare origine, quindi, a risultati esteticamente non gradevoli.[14] Si evita di fare la gemmazione sulle facce triangolari perché questa produrrebbe delle "valli" e non si otterrebbe un solido convesso.

---

[13] Sul sito https://mathcurve.com/polyedres/geode/geode.shtml è mostrata una particolare tassellazione delle facce di un icosaedro che porta alla costruzione della stessa cupola geodetica.

Anche se al tempo non le è stato attribuito questo aggettivo, viene spesso considerata come prima cupola "geodetica" quella (non più esistente) progettata dall'ingegnere ***Walther Bauersfeld*** nel 1922 e montata sulla sommità di uno degli edifici delle fabbriche di ottica Zeiss a Jena (Germania). Consisteva in una struttura semisferica a maglie triangolari (vedi Figura 13), rivestita di cemento spruzzato che doveva servire come superficie per testare il nuovo proiettore dello *Zeiss Planetarium*. Inaugurata e chiusa nel corso del 1924 [Ga1], viene successivamente smantellata e non ne restano che poche foto. Nel 1926 viene inaugurato il secondo planetario di Jena, progettato sempre da Bauersfeld; tale cupola ha un diametro di 23 metri e contiene una sala di 261 posti. Si tratta del più antico planetario ininterrottamente operativo al mondo. Bauersfeld progetterà altre cupole ma non di tipo geodetico[15].

Nel 1928 viene inaugurato il primo planetario di Roma, posto all'interno delle Terme di Diocleziano (ora sede del Museo Nazionale Romano), come riportato da A. Carlini e L. Tedeschini Lalli in [C-TL]. Lo schermo ed il proiettore Zeiss del planetario fanno parte del risarcimento danni della prima guerra mondiale. La griglia metallica dello schermo è tutt'ora visibile ed in [C-TL] viene evidenziata la similarità della triangolazione di questa cupola con quella (confrontabile solo in foto) della prima cupola di Jena; in particolare, si ha in entrambe l'inusuale presenza di vertici di grado 7, oltre a quelli di grado 5 e 6. Per questo motivo tali cupole non possono essere considerate geodetiche secondo la definizione data precedentemente e non possono essere ottenute proiettando una triangolazione di un solido platonico (o semiregolare) sulla sfera ad esso circoscritta. Anche O. Gaspar [Ga1, Ga2], che ha analizzato appunti e disegni di Bauersfeld recentemente ritrovati, mostra che sovrapponendo ad una foto della prima cupola di Jena le griglie ottenute suddividendo un icosaedro con 4 metodi diversi descritti in letteratura, la griglia che si adatta meglio corrisponde al metodo usato da Bauersfeld (così come ricostruito dai documenti ritrovati) ma non a quello indicato da Buckminster Fuller.

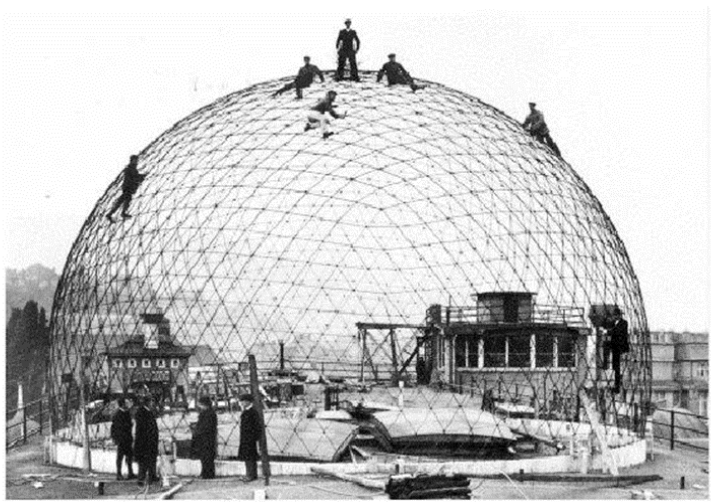

*Figura 13. La cupola di Jena in costruzione*
*https://en.m.wikipedia.org/wiki/File:Dome_Jena_UnderConstruction.jpg*

---

[14] Nel caso del dodecaedro, posto $R = 1$ il raggio della sfera circoscritta al poliedro, il suo spigolo è $\frac{\sqrt{5}-1}{\sqrt{3}}$, mentre il raggio della sfera inscritta è $\sqrt{\frac{5+2\sqrt{5}}{15}}$. Il lato obliquo della piramide risulta, quindi: $\sqrt{2 - 2\sqrt{\frac{5+2\sqrt{5}}{15}}}$ ed il rapporto fra questo lato e lo spigolo del dodecaedro vale circa 0,898.

[15] Poiché la prima cupola serve per un planetario, gli studi di Bauersfeld puntano ad una superficie il più possibile "liscia" e "sferica" da realizzare per mezzo di triangoli con aree il più possibile uguali e con poche differenze fra le lunghezze dei lati [Ga1, Ga2].

Alle cupole geodetiche viene usualmente associato il nome dal già citato architetto statunitense **Richard Buckminster Fuller** (1895 – 1983) perché è lui il primo a sviluppare in maniera sistematica l'idea di cupola geodetica a partire dal 1948 e a costruire i primi edifici con questa forma, ottenendone vari brevetti americani. Nelle estati del 1948 e 1949, Buckminster Fuller è "visiting professor" al Black Mountain College e lì con i suoi studenti sperimenta il primo esempio di cupola geodetica: una griglia semisferica di 48 piedi (circa 14,6 mt) di diametro, costruita con stecche di veneziane piegate in modo da seguire la curvatura delle geodetiche; questa struttura, però, collassa sotto il suo peso e viene ribattezzata *Supine Dome*. Nei tentativi successivi verranno usate aste rigide e varie soluzioni per le connessioni, ottenendo modelli più performanti. Nel 1951, Buckminster Fuller richiede il primo brevetto (concesso nel 1954) e sempre agli inizi degli anni '50 l'esercito americano gli commissiona varie cupole geodetiche da utilizzare sia come abitazioni sia come protezioni per strumentazioni costose. Nel 1954 viene invitato ad esporre alla Triennale di Milano, ottenendo così il primo riconoscimento pubblico a livello internazionale e nel 1967 raggiunge il massimo della popolarità realizzando la cupola per l'Expo di Montreal (Figura 9), di cui vedremo in seguito le caratteristiche geometriche.

Buckminster Fuller è stato un architetto autodidatta e "non convenzionale", un uomo carismatico dalla personalità eclettica, interessato a vari campi della Scienza, ma soprattutto inventore: gli vengono attribuiti ben 25 brevetti, anche se alcuni studiosi mettono oggi in discussione la piena paternità di alcuni risultati. I suoi studi e le sue invenzioni sono tutti legati alla ricerca e alla successiva messa in pratica di principi strutturali che governano l'Universo. Come spiega nella raccolta *Synergetics* [BF], tutte le strutture, interpretate in maniera appropriata, dal sistema solare all'atomo, sono strutture di **tensegrità**. Con il termine da lui coniato *tensegrity*, Buckminster Fuller contrae le parole *tension* e *integrity* per descrivere il principio che governa tutto l'Universo di *compressione isolata* e *tensione continua:* componenti isolate che "lavorano" in compressione (di solito barre o puntoni nelle realizzazioni) sono all'interno di un sistema continuo in tensione (in genere realizzato con cavi o tendini). Le componenti in compressione non si toccano l'una con l'altra; lo scultore Kenneth Snelson, allievo di Buckminster Fuller e primo[16] a realizzare e poi sviluppare *Tensegrity Sculptures*, descrive le realizzazioni come barre che galleggiano nell'aria senza alcun supporto "solido". Le strutture devono, quindi, la loro stabilità ad un equilibrio di forze opposte che si compensano e Buckminster Fuller individua in questo meccanismo un principio fondamentale della Natura che, con pochi elementi, riesce a realizzare una struttura "forte". Sempre prendendo spunto dalla Natura, Buckminster Fuller definisce una geodetica come la relazione più economica fra due eventi e dichiara che la geometria delle cupole geodetiche è la più efficiente fra tutti i sistemi a cupola triangolati a curvatura composta. Afferma anche che le cupole geodetiche sono tutte strutture di tensegrità, indipendentemente dalla visibilità per l'osservatore della differenziazione fra tensione e compressione[17]. Poiché queste strutture utilizzano una quantità minima di materiale, sono secondo Buckminster Fuller le più leggere ed economiche.

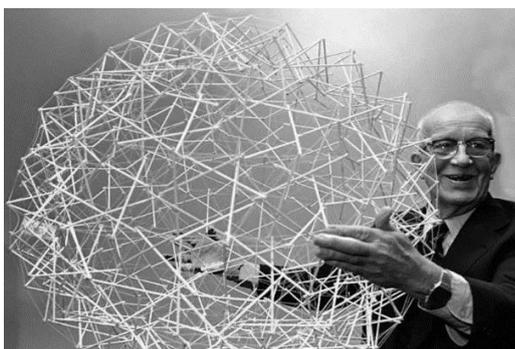

*Figura 14. R. Buckminster Fuller con una struttura Tensegrity Sphere - 1979*[18]

---

[16] Tre persone hanno un brevetto che riguarda strutture di tensegrità: R. Buckminster Fuller (1962, anche se lui stesso ammette di non averle mai costruite), G. Emmerich (1964, artista ungherese che sembra aver lavorato senza conoscere Buckminster Fuller e Snelson), Kenneth Snelson (1965), che Buckminster Fuller inizialmente riconosce come primo realizzatore, ma poi sembra "dimenticarsene".

[17] In realtà le cupole realizzate da Buckminster Fuller, come quella di Montreal, non sembrano essere strutture di tensegrità.

[18] **https://www.flickr.com/photos/poetarchitecture/26806590126/in/photolist-GQNMjo-hESW2z-GMT4BP-ejcfv3-criycW-r4RXrm-qixJV2-3ZnJR-3ZnKg-5mMEfE-5mHpSD-5mMEDd-VR9y-VR7Y-VR9e-VR7D-**

La facile trasportabilità di queste cupole le rende interessanti anche per l'esercito americano, che, come accennato prima, ne installa centinaia in Canada per proteggere postazioni radar. Negli anni '60 le conferenze di Buckminster Fuller sono molto seguite dalle giovani generazioni americane e vengono costruite nuove comunità (come *Drop City* in Colorado), animate anche dalle idee di organizzazione sociale dell'architetto [Tu]; tuttavia, già un decennio dopo queste comunità sono sparite quasi tutte. A differenza di quanto ipotizzato da Buckminster Fuller, le cupole geodetiche non risolvono i gravi problemi abitativi degli Stati Uniti di quegli anni; comunque, sono ancora fonte di ispirazione e ricerca per architetti e ingegneri.

Il principio di tensegrità è ad oggi il modello di riferimento nella biomeccanica. A partire dagli anni '80, il biologo Donald Ingber, studiando la biomeccanica dei sistemi viventi, mostra che la Natura usa il principio di tensegrità per stabilizzare la forma delle cellule viventi e determinarne la risposta alle sollecitazioni meccaniche. Le strutture di tensegrità, scrive Ingber, sono efficienti e flessibili mantenendo la robustezza [In].

## 4 - Dalle cupole geodetiche ai virus passando per la PopArt

Una diversa "triangolazione" lega le cupole geodetiche con l'arte e con la biologia, come dettagliato da Gregory J. Morgan in diversi lavori [Mo1, Mo2, Mo3].

Nel 1956 Francis Crick e James Watson [C-W] pubblicano su *Nature* un famoso articolo in cui ipotizzano una struttura poliedrica del capside, cioè del rivestimento proteico esterno che racchiude l'RNA, di molti virus comunemente detti *sferici*. Come già suggerito precedentemente da altri scienziati, questi virus presentano simmetrie che sono possibili in tetraedri, ottaedri e icosaedri. Poiché il genoma virale è "povero", Crick e Watson ipotizzano che le (poche) proteine del capside siano organizzate in sottounità asimmetriche posizionate nello stesso modo su ogni *faccia* di un poliedro; gli scienziati pensano all'icosaedro perché questo consente il posizionamento di 60 unità asimmetriche (contro le 24 dell'ottaedro e le 12 del tetraedro)[19] su ogni faccia ed anche perché, fra i poliedri considerati, l'icosaedro racchiude un volume più grande (e quindi una maggior quantità di acido nucleico). Questa teoria trova inizialmente scettici diversi esperti, anche di spicco, del settore[20] mentre altri scienziati, proseguendo lo studio della diffrazione ai raggi X di alcuni virus cristallizzati, ipotizzano l'esistenza di un numero anche maggiore di sottounità e richiedono, quindi, un nuovo modello. In particolare, i risultati di John Thomas Finch e Aaron Klug sul virus della polio vengono diffusi anche da alcuni quotidiani e l'artista John McHale se ne interessa. McHale è al tempo un membro dell'*Independent Group* di Londra, un gruppo influente di artisti "ribelli" in cerca di nuove fonti di ispirazione estetica da trarre dalla cultura popolare, dalla produzione di massa e dalla scienza. L'*Independent Group* è considerato l'antesignano in Inghilterra della *PopArt* (termine coniato secondo alcuni storici dallo stesso McHale ma attribuito dalla maggioranza a Lawrence Alloway, un altro esponente del gruppo). McHale è già in contatto con Buckminster Fuller, contribuendo alla diffusione delle sue tematiche e dei suoi progetti; leggendo gli articoli sul virus della polio, pensa che ci possa essere un legame con le cupole geodetiche ed organizza un incontro fra Buckminster Fuller e Klug. Anche se il primo incontro porta solo alla condivisione di simmetrie icosaedriche, Klug continua a visionare i disegni di cupole geodetiche presenti nei manoscritti di Buckminster Fuller e studia anche le *Tensegrity Structures* realizzate da Kenneth Snelson. In particolare, Klug è colpito dalla *270-struts tensegrity sphere*; formata da 270 puntoni e, quindi, 540 sottounità ottenute dividendo ogni puntone nelle sue due parti asimmetriche, questa scultura fornisce un modello per organizzare $540 = 60 \cdot 9$ unità proteiche asimmetriche su una sfera. Negli stessi anni il collega Donald L.D. Caspar, che ha già presentato dei modelli di capside a simmetria icosaedrica costruiti con palline da pingpong, trascorre

**VR8M-8y9tDo-8y6sNX-qnhPRv-sSPR3B-ta1L5A-sSFpTo-t7XFvh-t7Xf6u-t7WDZd-t7W8aY-sSFCyf-t7WNX3-sdgce7-sSGbAS-sSEAJd-sSH5eG-t7WeNY-sdsw7p-sdrtJa-t7WvQs-ta2Hj3-taiBsF-tagNuP-sSPTcM-t7WCsq-ta1wys-sSNNhP-ta2Tpo-sSFMmJ-sSPk8M-sdrEH4-ta2Jc5-sSHcrN**

[19] Dalle osservazioni sperimentali, gli scienziati ipotizzano che le sottounità proteiche (asimmetriche) siano aggregate rispettando le simmetrie di un gruppo puntuale cubico (gruppo delle simmetrie rotazionali di un cubo o di un ottaedro); per avere la stessa disposizione su ogni faccia del poliedro, conservandone le simmetrie, si hanno le limitazioni descritte. Si veda anche la Nota 25.

[20] I gruppi di simmetria dei cristalli non prevedono simmetrie di ordine 5 ed anche in Natura queste sono abbastanza rare; i *quasicristalli* verranno scoperti a metà degli anni Ottanta.

del tempo con Buckminster Fuller; identificando come possibili sottounità proteiche le facce di una cupola geodetica, determina la formula (che descriveremo nel paragrafo 5) $T = m^2 + mn + n^2$, con $m$ e $n$ interi non negativi e non entrambi nulli, che porta alla definizione di $T$ come numero di triangolazione o numero del virus: ogni virus *sferico* può essere classificato dal suo numero $T$ ed il relativo capside è formato da $60T$ sottounità, $3T$ per ogni sfaccettatura. Introducendo il concetto di *quasi-equivalenza* per i legami delle sottounità proteiche, Caspar e Klug pubblicano i loro risultati nel 1962 [C-K]; questo articolo, insieme a quello di Crick e Watson, è considerato una delle "pietre miliari" della nascente biologia molecolare. L'analogia del capside con le cupole geodetiche non si ha solo nella struttura *triangolata* di ogni faccia, ma anche nel fatto che i legami di *quasi-equivalenza* comportano che le sfaccettature non sono identiche fra loro, così come le facce triangolari di una cupola geodetica non sono identiche fra loro. Caspar posiziona le sottounità non sui vertici ma sulle facce triangolari di una cupola geodetica; quindi, di fatto lavora con il poliedro duale, ovvero un poliedro di Goldberg. Come scoprirà più tardi, la formula per calcolare $T$ è già presente in [Go2]. Per un confronto fra le varie costruzioni si veda anche [B-G-S] e per le analogie fra cupole geodetiche e struttura dei virus [Co].

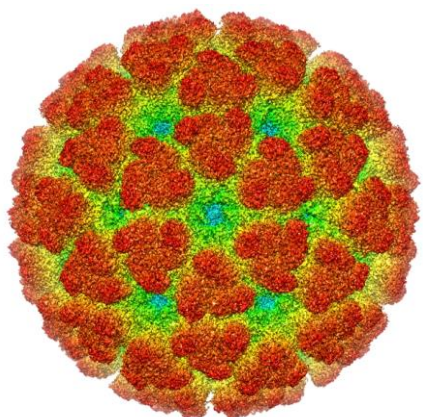

*Figura 15. Chikungunya Virus   T=4, (m,n)=(2,0)*
*Foto di A2-33 - Opera propria, CC BY-SA 3.0, https://commons.wikimedia.org/w/index.php?curid=30051546*

La scoperta, in anni successivi, di strutture virali a simmetria icosaedrica (come il virus della polio), che non rientrano nel modello di Caspar e Klug, ha portato alla recente *Viral Tiling Theory*, in cui la struttura del capside viene descritta con tassellazioni più generali (tipo Penrose) attraverso i gruppi matematici dei *quasicristalli*. Le ricerche in tal senso sono iniziate con Reidun Twarock [Tw] e proseguono tutt'ora.

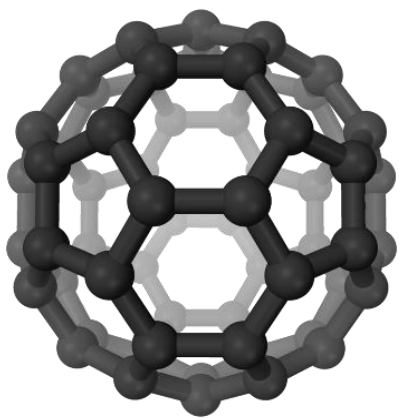

*Figura 16.  Buckminsterfullerene*
*Foto di Benjah-bmm27 - Own work, Public Domain https://commons.wikimedia.org/w/index.php?curid=1913689*

Nel 1985 un gruppo di 5 chimici inglesi e statunitensi pubblica un articolo [K-H-OB-C-S] intitolato *$C_{60}$: Buckminsterfullerene*, in cui viene descritta una nuova molecola costituita esclusivamente da atomi di carbonio. Si tratta di un aggregato particolarmente stabile formato da 60 atomi di carbonio, nato in laboratorio a seguito di esperimenti sulla frammentazione della grafite in "condizioni violente" che simulano quelle che accadono nell'Universo. Gli autori dell'articolo spiegano che, dovendo ipotizzare una struttura per questa specie e considerando che il carbonio non può avere un intorno tetraedrico come nel diamante, hanno pensato ad una forma "sferica" e si sono così rivolti agli studi di Buckminster Fuller. Fra le figure con alta simmetria, l'ipotesi dell'icosaedro troncato, con 60 atomi di carbonio disposti negli altrettanti vertici, è stata ritenuta "di bellezza unica". A seguito di queste scoperte, nel 1996 Harold Walter Kroto, Robert Curl e Richard Smalley ricevono il Premio Nobel per la Chimica. Sempre nell'articolo citato, gli autori precisano che una struttura così importante merita un proprio nome e, dopo aver considerato alternative quali "soccerene" o "sferene", scelgono *Buckminsterfullerene*, talvolta abbreviato in *Buckyball*. Col termine *fullereni* vengono oggi indicati gli elementi appartenenti a una classe di allotropi del carbonio diversi dalla grafite e dal diamante (gli allotropi sono modificazioni strutturali diverse dello stesso elemento), i cui atomi sono tri-coordinati e formano strutture a gabbia chiusa; il *Buckminsterfullerene* è il più stabile fra questi. Studi recenti fatti sulle simmetrie di capside virale indicano che gli stessi modelli geometrici possono essere applicati alle *carbon onions*, nanoparticelle di solo carbonio costituite da gusci concentrici di molecole di fullerene.

## 5 - *Primo tipo di triangolazione di un icosaedro*

Descriviamo ora uno dei metodi più comuni per costruire le sfere geodetiche in modo che queste abbiano la maggior regolarità possibile, partendo da un icosaedro (vedi anche [C-P-T1]).

Dividiamo i lati di ogni sua faccia in $m$ segmenti congruenti, con $m$ numero naturale, e uniamo i punti così ottenuti con segmenti paralleli ai lati di ciascuna faccia. Quindi, ogni faccia risulta tassellata con $m^2$ triangoli equilateri congruenti fra loro, come si può vedere nelle figure seguenti (se $m = 1$ non si ha alcuna tassellazione); si dice che la suddivisione ha **frequenza** $m$v oppure che è di **tipo** $m$v.

Questo modo di tassellare ciascuna faccia dell'icosaedro con triangoli equilateri è detto anche **ripartizione alternata**[21].

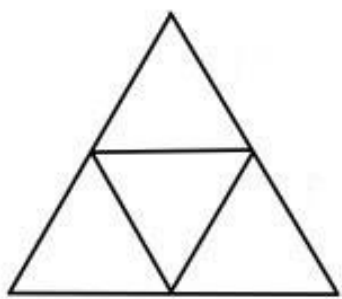
*Figura 17. Frequenza* 2v

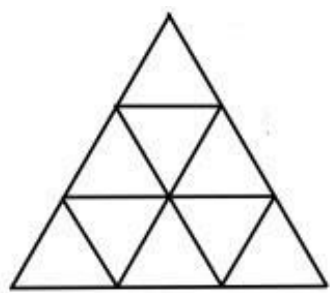
*Figura 18. Frequenza* 3v

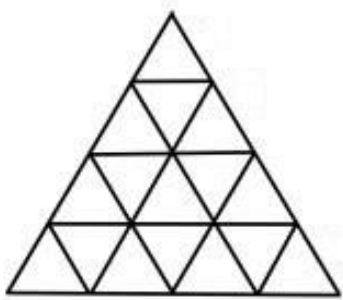
*Figura 19. Frequenza* 4v

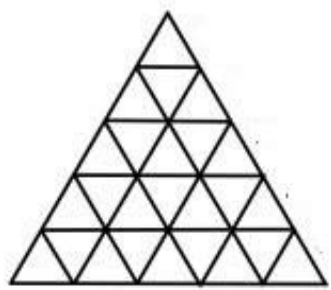
*Figura 20. Frequenza* 5v

Dal centro della sfera circoscritta all'icosaedro, proiettiamo i vertici dei triangoli così ottenuti sulla sfera stessa; otteniamo in tal modo una sfera geodetica a facce triangolari.

Ad esempio, la sfera geodetica 2v (nella Figura 21) ha 80 facce, 42 vertici (12 di grado 5 e 30 di grado 6) e 120 spigoli.

---

[21] Storicamente, questo tipo di ripartizione nasce dopo quella oggi detta di tipo ($m$, 0) (vedi il paragrafo successivo) e per questo viene indicata come "alternativa" [Ke].

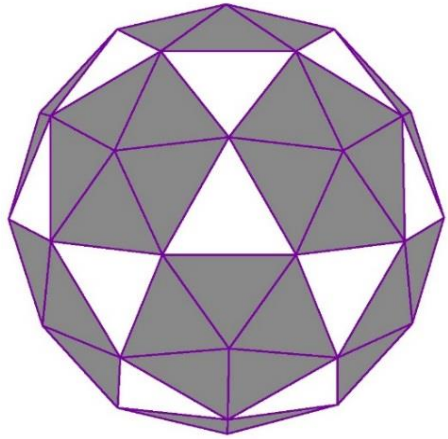
*Figura 21. Sfera geodetica 2v*

Il poliedro di Goldberg, duale di questa sfera geodetica, è il **triacontaedro rombico troncato**, che ha 80 vertici, 120 spigoli e 42 facce: 30 esagonali e 12 pentagonali. Poiché questo poliedro ha una forma simile all'icosaedro troncato, esso è anche chiamato *super pallone da calcio.*

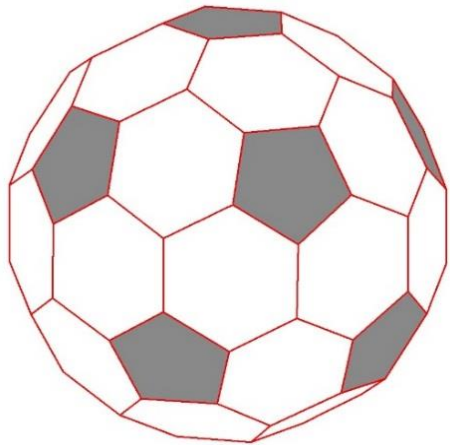
*Figura 22. Triacontaedro rombico troncato*

Per una lista di sfere geodetiche e dei loro duali rimandiamo a Magnus J. Wenninger [We].

Una costruzione analoga a questa può essere ottenuta anche partendo da altri poliedri regolari; tuttavia la maggior parte delle sfere geodetiche sono costruite partendo dall'icosaedro perché questo è il poliedro regolare che, a parità di volume, ha una superficie più piccola.

Operando una suddivisione di tipo $m$v sull'icosaedro, si ottengono $20m^2$ triangoli equilateri più piccoli e quindi una sfera geodetica con $20m^2$ facce triangolari. Anche se ogni faccia dell'icosaedro è tassellata con triangoli equilateri, le facce della sfera geodetica ottenuta sono triangoli isosceli (solo alcuni equilateri) non tutti uguali fra loro. Ad esempio, nella sfera geodetica 3v compaiono tre diverse lunghezze di spigoli, in quella 5v nove lunghezze diverse. Nel corso degli anni sono state sviluppate varie tecniche di calcolo e, recentemente, anche software *ad hoc.*

Per calcolare la frequenza di una sfera geodetica (ottenuta a partire da un icosaedro) si considerano due vertici "vicini" di grado 5 (che corrispondono agli estremi di uno spigolo dell'icosaedro) e si contano i segmenti che uniscono tali vertici: il numero di questi segmenti corrisponde alla frequenza.

Nei seguenti disegni sono presentate alcune *semisfere* geodetiche con frequenze diverse; in questi casi le facce triangolari sono $10m^2$.

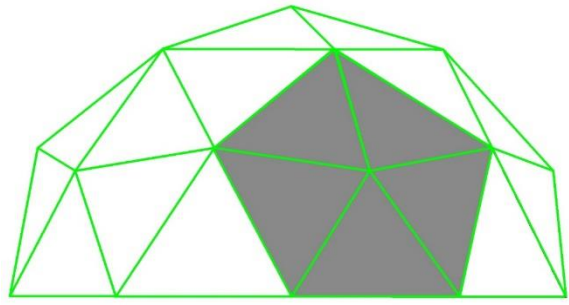
*Figura 23. Semisfera geodetica 2v*

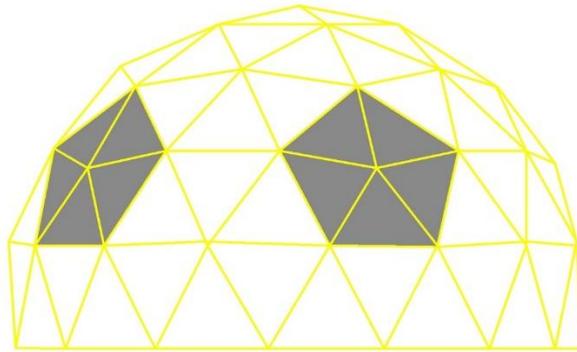
*Figura 24. Semisfera geodetica 3v*

Notiamo che la sfera geodetica 3v corrisponde a quella ottenuta dall'innalzamento dell'icosaedro troncato (vedi la Figura 12).

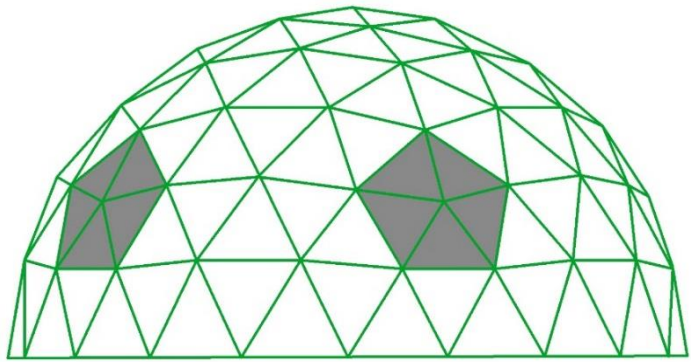
*Figura 25. Semisfera geodetica 4v*

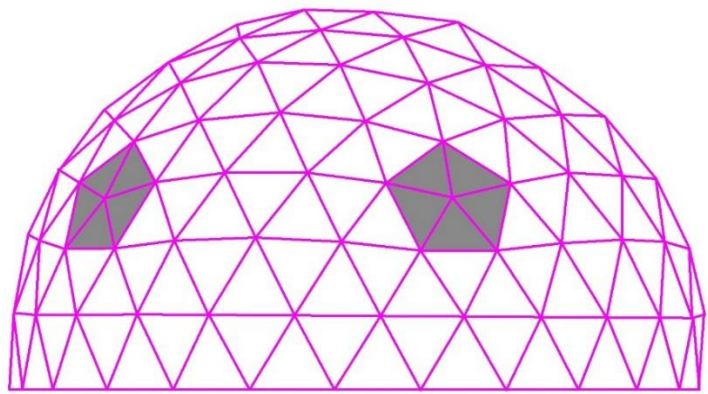
*Figura 26. Semisfera geodetica 5v*

Nel Parco de la Villette a Parigi si trova il **Géode**, una grande sfera (quasi completa) geodetica d'acciaio di 36 metri di diametro, realizzata nel 1985 dall'architetto Adrien Fainsilbe. Al suo interno si trova una sala cinematografica di 400 posti, dotata di uno schermo semisferico di 1000mq, il più grande al mondo. La parte esterna della sfera ha un'altezza di circa 30 metri. Si tratta di una sfera geodetica di tipo 20v con la parte esterna formata da circa 6433 triangoli di acciaio inossidabile.

La **Biosphère** a Montreal (Canada), menzionata precedentemente, è la cupola geodetica[22] che Richard Buckminstrer Fuller ha progettato come Padiglione degli Stati Uniti per l'Esposizione Universale del 1967. La cupola è di tipo 16v ed è circa i 4/5 di una sfera di diametro 76,2 metri, perciò è alta circa 60,96 metri ed è formata da circa 4096 triangoli [B-B].

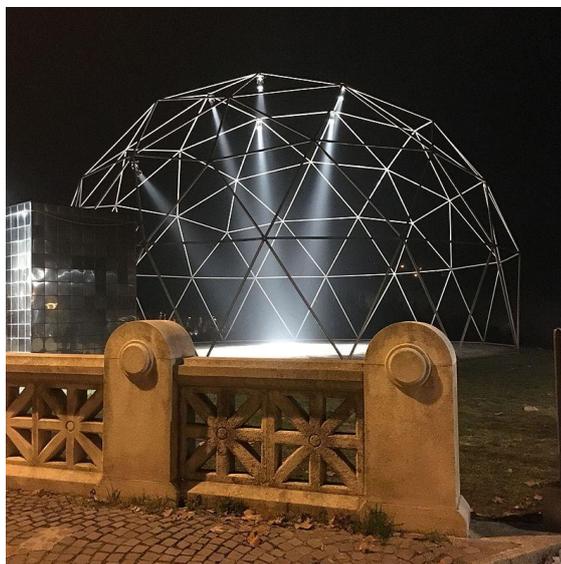

*Figura 27. Semisfera geodetica di Buckminster Fuller a Spoleto  Foto di Manuela Rosi [licenza CC]*
*https://commons.wikimedia.org/w/index.php?curid=45968494*

Nel 1967 Buckminstrer Fuller ha donato alla città di Spoleto, in occasione del *X Festival dei Due Mondi*, una semisfera geodetica chiamata **Spoletosfera.** Essa è formata da elementi metallici, misura 21 metri di diametro, è di tipo 3v ed è formata da 90 facce triangolari. Con opportune modifiche (ad esempio, coprendola con un telo) tale struttura può essere utilizzata per organizzare conferenze, mostre, spettacoli.

## 6 - Secondo tipo di triangolazione di un icosaedro

Descriviamo un altro metodo, introdotto da Michel Goldberg in [Go2] e ripreso da Donald L. D. Caspar e Aaron Klug [C-K], per realizzare una famiglia infinita di sfere geodetiche che possiedono le simmetrie rotazionali dell'icosaedro e dotate della massima regolarità possibile (vedi anche [C-P-T2]).

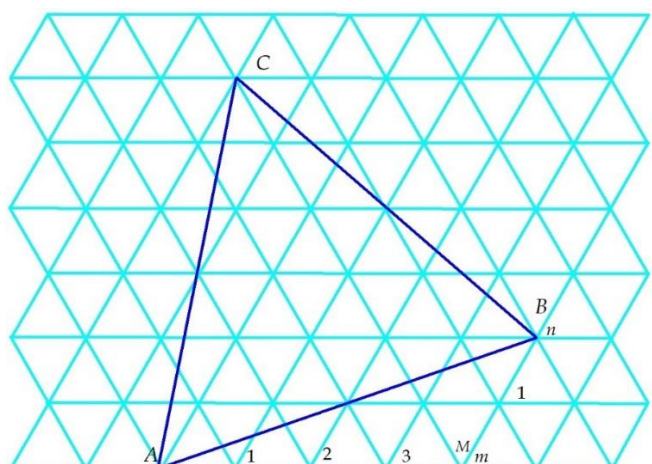

*Figura 28. La suddivisione m = 4, n = 2*

---

[22] Sono stati due ex studenti di Buckminster Fuller, Jeffrey Lindsay e Don Richter, a progettare nel 1950 il prototipo della cupola realizzata a Montreal, usando il metodo che hanno chiamato "della griglia regolare".

Consideriamo una tassellazione del piano formata da triangoli equilateri di lato unitario, come in Figura 28. Consideriamo due numeri naturali $m$ e $n$ non entrambi nulli. A partire da un vertice $A$ della tassellazione, percorriamo $m$ spigoli nella direzione "orizzontale" del reticolo e poi $n$ spigoli nella direzione ruotata, rispetto alla precedente, di 60° in senso antiorario per arrivare al vertice $B$; nella Figura 28, $m = 4$ e $n = 2$. Costruiamo il triangolo equilatero $ABC$; notiamo che il vertice $C$ è un vertice della tassellazione e può essere ottenuto, a partire da B, percorrendo $m$ spigoli nella direzione ruotata di 60° rispetto all'ultima percorsa e poi $n$ spigoli nella direzione ulteriormente ruotata di 60° (sempre in senso antiorario). Ripetendo lo stesso procedimento a partire dal vertice $C$ si ritrova il punto $A$.

La costruzione precedente è detta tassellazione di **tipo** ($m$, $n$) del triangolo $ABC$.

Notiamo che, eccetto i casi $m = 1$ e $n = 0$ oppure $m = 0$ e $n = 1$, il triangolo equilatero $ABC$ è più grande di quelli del reticolo di partenza.

Sostituiamo ad ogni faccia di un icosaedro questo triangolo $ABC$ con tutti gli elementi della tassellazione originale in esso contenuti. Successivamente, dal centro della sfera circoscritta all'icosaedro proiettiamo i **vertici** dei triangoli della tassellazione sulla sfera stessa; vengono proiettati solo i vertici dei triangoli e non eventuali altri punti che si trovano sugli spigoli del poliedro, eccetto i casi $m = 0$ oppure $n = 0$. Si ottiene così una sfera geodetica $S$ a facce triangolari. Notiamo che i vertici dei triangoli piccoli possono trovarsi in parte su una faccia dell'icosaedro e in parte su una faccia contigua; tuttavia, i triangoli che si trovano in parte su una faccia e in parte su un'altra faccia contigua, si "ricompongono" nella proiezione.

I triangoli della tassellazione una volta proiettati non sono più equilateri[23] in quanto i loro vertici, prima della proiezione, non appartengono alla sfera circoscritta all'icosaedro a meno che non siano vertici dell'icosaedro stesso. Le facce triangolari della sfera geodetica ottenuta con questa costruzione sono comunque triangoli isosceli; inoltre, esse non sono tutte uguali fra loro.

Notiamo che tutti i vertici dell'icosaedro sono anche vertici della tassellazione, quindi appartengono alla sfera geodetica. Intorno a questi vertici, si dispongono cinque triangoli piccoli, mentre intorno a tutti gli altri vertici della sfera geodetica si collocano sei triangoli piccoli.

La sfera geodetica così ottenuta è detta di **tipo** ($m$, $n$) perché ottenuta con una tassellazione di tipo ($m$, $n$) delle facce dell'icosaedro.

Notiamo che le tassellazioni di tipo ($m$, 0) o di tipo (0, $n$) sono formate da triangoli equilateri con i lati paralleli ai lati del triangolo equilatero $ABC$ (vedi la Figura 29). In questi casi (e solo in questi!) tutti i punti che vengono a trovarsi sugli spigoli dell'icosaedro corrispondono a vertici dei triangoli della tassellazione e vengono quindi proiettati sulla sfera. Le sfere geodetiche di tipo ($m$, 0) o (0, $n$) sono dette di **Classe I** e corrispondono a quelle viste nella costruzione illustrata nel paragrafo precedente. Nella figura seguente è riportata una tassellazione di tipo (5, 0).

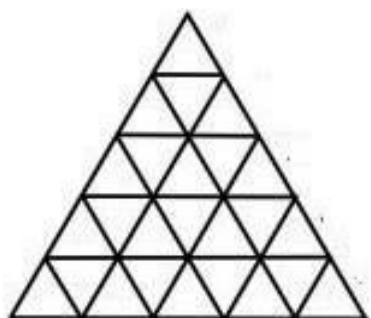

*Figura 29. Tassellazione di tipo (5, 0)*

La semisfera geodetica corrispondente alla precedente tassellazione è rappresentata nella Figura 30.

---

[23] Notiamo che, oltre ai tre poliedri regolari a facce triangolari (tetraedro, ottaedro, icosaedro), esistono soltanto 5 poliedri convessi le cui facce sono triangoli equilateri; tali poliedri sono chiamati **deltaedri**. I deltaedri sono stati studiati dal matematico olandese Bartel Leendert van der Waerden nel 1947; sono stati chiamati così da Martyn Cundy nel 1952 [Cu]. Il loro nome deriva dalla lettera greca *Delta*, che ha la forma di un triangolo.

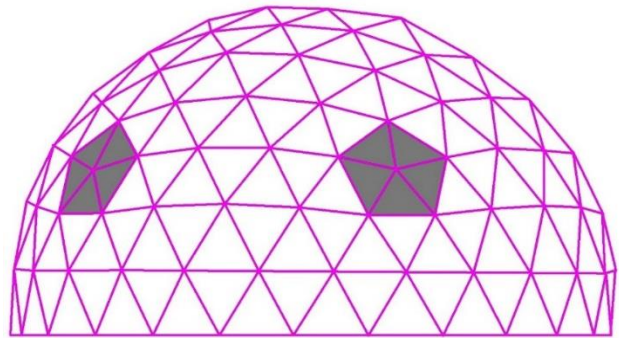
*Figura 30. Semisfera geodetica di tipo (5, 0)*

Nel caso di tassellazioni di tipo (*n*, *n*), con $n \neq 0$, si dice che la cupola geodetica è di **Classe II**. Nella Figura 31 è disegnata una tassellazione di tipo (2, 2).

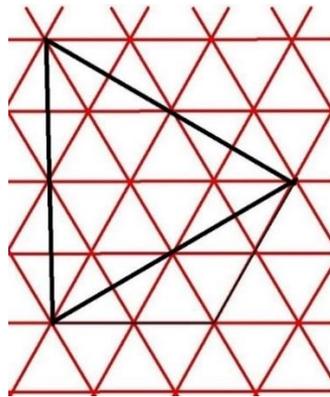
*Figura 31. Tassellazione di tipo (2, 2)*

La sfera geodetica corrispondente alla precedente tassellazione è rappresentata nella Figura 32.

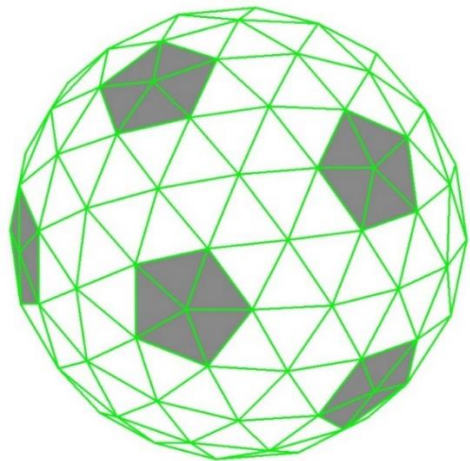
*Figura 32. Sfera geodetica di tipo (2, 2)*

Nel sito internet https://mathcurve.com/polyedres/geode/geode.shtml è riportato un esempio di cupola geodetica del tipo (2, 2), usata per telecomunicazioni, che si trova nell'Antartide.

È interessante notare che il pentacisdodecaedro inscritto in una sfera non è altro che la cupola geodetica di tipo (1, 1).

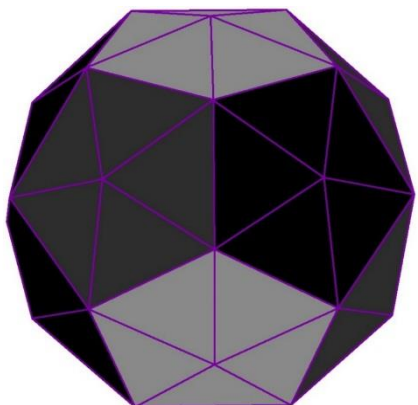
*Figura 33. Pentacisdodecaedro*

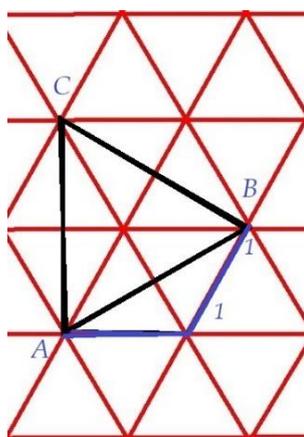
*Figura 34. Tassellazione di tipo (1, 1)*

Come si può vedere anche dalle Figure 31 e 34, le tassellazioni di tipo di tipo ($n$, $n$) sono formate da triangoli equilateri e *semitriangoli* equilateri con ciascun lato perpendicolare ad uno dei lati del triangolo *ABC*.

Nel caso di tassellazioni di tipo ($m$, $n$), con $m \neq n$ e $m, n \neq 0$, si dice che la cupola geodetica è di **Classe III**. Notiamo che i poliedri di questa classe sono chirali.

Ogni tassellazione di questo tipo è formata da triangoli equilateri, o porzioni di triangoli equilateri, i cui lati non sono né paralleli né perpendicolari ai lati del triangolo *ABC,* come si vede, ad esempio, considerando la tassellazione di tipo (2, 1).

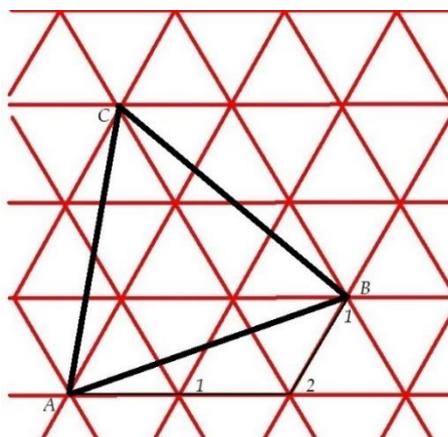
*Figura 35. La tassellazione di tipo (2, 1)*

Con questa tassellazione otteniamo la sfera geodetica della Figura 36.

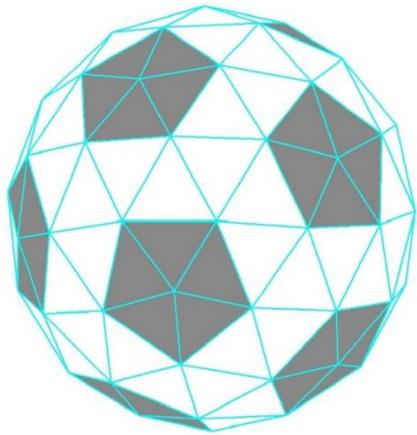

*Figura 36. La sfera geodetica di tipo (2, 1)*

La costruzione fatta nel caso *m* = 0 o *n* = 0 può essere eseguita anche a partire dagli altri poliedri regolari a facce triangolari, cioè il tetraedro e l'ottaedro ([Ba], pp. 221-223); tuttavia, quasi tutte le sfere geodetiche si ottengono considerando l'icosaedro come poliedro di partenza.

Calcoliamo ora il numero di facce triangolari presenti in una cupola geodetica.

Facciamo riferimento alla Figura 28 e ricordiamo che i triangoli equilateri che formano la tassellazione del piano hanno il lato unitario.

Determiniamo anzitutto il modulo del vettore $\overrightarrow{AB}$. Si ha $\overrightarrow{AB} = \overrightarrow{AM} + \overrightarrow{MB}$. Supponiamo di riferire il piano ad un sistema di coordinate cartesiane ortogonali, in cui il punto *A* è l'origine e la retta *AM* l'asse delle ascisse, con verso concorde con quello del vettore $\overrightarrow{AM}$.

Si ottiene $\overrightarrow{AM} = (m, 0)$ e $\overrightarrow{MB} = \left(\frac{1}{2}n, \frac{\sqrt{3}}{2}n\right)$. Ne segue che

$$\overrightarrow{AB} = \overrightarrow{AM} + \overrightarrow{MB} = \left(m + \frac{1}{2}n, \frac{\sqrt{3}}{2}n\right)$$

quindi la misura del segmento *AB* è

$$\sqrt{\left(m + \frac{1}{2}n\right)^2 + \frac{3}{4}n^2} = \sqrt{m^2 + mn + n^2}$$

Ne consegue che l'area del triangolo equilatero *ABC* è: $\frac{\sqrt{3}}{4}(m^2 + mn + n^2)$.

Il numero $T$ di triangoli della tassellazione presenti in *ABC* si ottiene dividendo l'area di *ABC* per l'area di un triangolino della tassellazione, che è $\frac{\sqrt{3}}{4}$, quindi

$$T = m^2 + mn + n^2.$$

Di conseguenza il numero delle facce triangolari che formano una cupola geodetica di tipo (*m*, *n*) è 20$T$.

Il numero $T$ è anche chiamato **numero di triangolazione** e, come accennato nel paragrafo 4, è molto utilizzato in virologia. Ad esempio, il numero delle facce della sfera geodetica nella Figura 36 è 140. Invece il numero delle facce della sfera geodetica di tipo (1, 1) (cioè il pentacisdodecaedro) è 60.

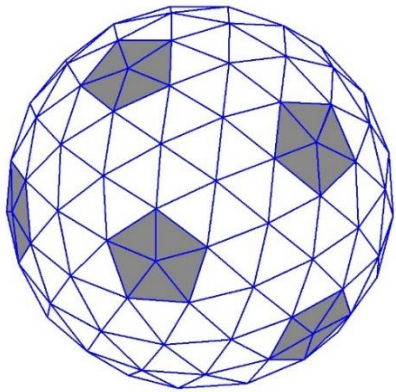

*Figura 37. Sfera geodetica di tipo (3, 1)*

La sfera geodetica di tipo (3,1) (vedi Figura 37) ha 260 facce.
Nel paragrafo successivo calcoleremo anche il numero degli spigoli e dei vertici delle sfere geodetiche.

## 7– Ulteriori considerazioni sui poliedri di Goldberg

Come abbiamo già osservato, i poliedri duali di una cupola geodetica sono poliedri di Goldberg. Ad esempio, il duale della cupola geodetica del tipo (1, 1) è l'**icosaedro troncato**.

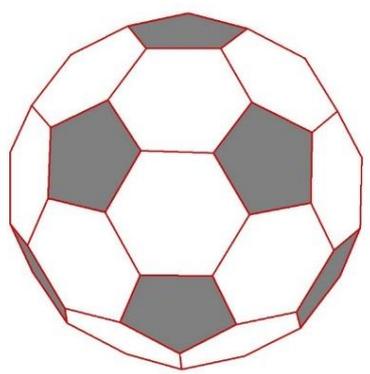

*Figura 38. Icosaedro troncato (pallone da calcio)*

In alcuni testi, anche i poliedri di Goldberg sono considerati sfere geodetiche; questo ci sembra poco corretto perché, a parte il dodecaedro regolare e l'icosaedro troncato, i poliedri di Goldberg non sono inscrivibili in una sfera.

Esistono comunque diverse strutture architettoniche che hanno la forma di poliedri di Goldberg. Una delle più interessanti è l'*Eden Project*, un complesso turistico in Cornovaglia, a circa 2 km dalla città di St Austell, ricavato nello spazio di una ex-cava. In questo complesso si trovano due grandi biosfere il cui interno ospita circa 100.000 piante provenienti da tutto il mondo.

Come abbiamo detto, gli spigoli dei poliedri di Goldberg sono tutti uguali, ma i poligoni che formano le sue facce esagonali non sono regolari, a parte l'icosaedro troncato.

Abbiamo già osservato che nei poliedri di Goldberg si ha $F_5 = 12$, da cui

$$F_6 = \frac{V}{2} - 10$$

Nel seguito indicheremo con il simbolo $GC_{m,n}$ il poliedro di Goldberg duale della sfera geodetica del tipo $(m, n)$.

Il poliedro duale della cupola geodetica del tipo (2, 1) della Figura 36 è il poliedro di Goldberg $GC_{2,1}$ della Figura 39.

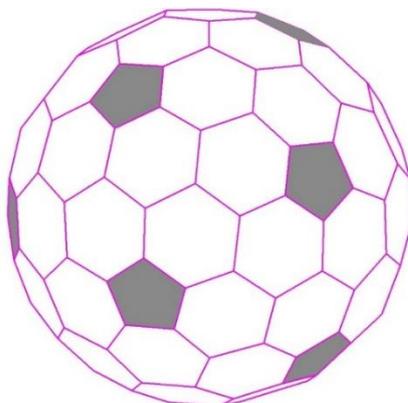

*Figura 39. Poliedro di Goldberg GC$_{2,1}$*

Abbiamo visto che la sfera geodetica di tipo (2, 1) ha **T** = 7 e quindi 20**T** = 140 facce; di conseguenza il poliedro di Goldberg $GC_{2,1}$ ha 140 vertici. Calcoliamo gli altri elementi di questo poliedro.

Dalla formula $F_6 = V/2 - 10$, si ha $F_6 = 60$; sommando a queste le 12 facce pentagonali, si ha che $GC_{2,1}$ ha 72 facce; quindi, dalla relazione di Eulero segue che ha 210 spigoli. Allora la sfera geodetica di tipo (2, 1), duale del poliedro di Goldberg $GC_{2,1}$, ha 210 spigoli e 72 vertici.

Generalizzando questo procedimento, si ottengono queste formule per le sfere geodetiche di tipo (*m*, *n*):

$$F = 20\boldsymbol{T} \qquad S = 30\boldsymbol{T} \qquad V = 10\boldsymbol{T} + 2$$

dove, ricordiamo, $\boldsymbol{T} = m^2 + mn + n^2$.

La prima relazione è già stata dimostrata nel paragrafo precedente. Per provare la seconda, ricordiamo che nei poliedri di Goldberg il numero degli spigoli è $S = \frac{5F_5 + 6F_6}{2}$ e che tale numero è uguale al numero degli spigoli della sfera geodetica sua duale; inoltre, il numero dei vertici $V_G$ di un poliedro di Goldberg soddisfa la seguente relazione

$$F_6 = \frac{V_G}{2} - 10$$

dove $V_G$ = numero delle facce della sfera geodetica duale = 20**T**. Dunque, dalla relazione precedente si ha:

$$S = \frac{5 \cdot 12}{2} + \frac{6 \cdot F_6}{2} = 30 + 3\left(\frac{V_G}{2} - 10\right) = \frac{3}{2} \cdot 20\boldsymbol{T} = 30\boldsymbol{T}$$

Dalla relazione di Eulero segue la terza formula: $V = 10\boldsymbol{T} + 2$.

Ad esempio, la sfera geodetica di tipo (4, 1) (vedere la Figura 40) ha **T** = 21. Di conseguenza essa ha 420 facce, 630 spigoli e 212 vertici. Segue che il corrispondente poliedro di Goldberg $GC_{4,1}$ (vedere la Figura 41) ha 420 vertici, 630 spigoli e 212 facce.

È interessante notare il fatto che nelle sfere geodetiche e nei poliedri di Goldberg è sufficiente conoscere uno soltanto fra i numeri *S, F,* e *V* per determinare gli altri due; invece, in un qualunque poliedro, occorre conoscerne due per determinare il terzo, usando la relazione di Eulero.

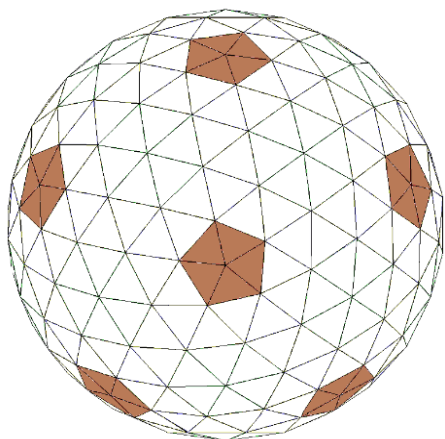

*Figura 40. Sfera geodetica di tipo (4, 1)*

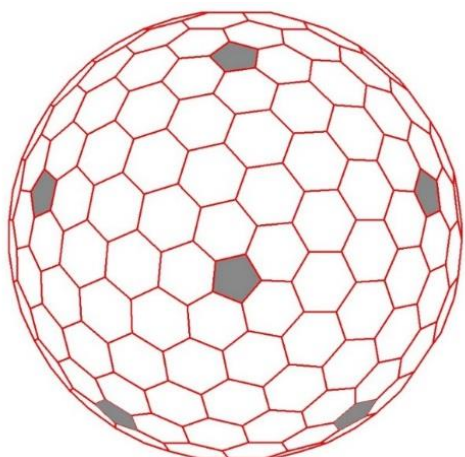

*Figura 41. Poliedro di Goldberg GC$_{4,1}$*

**Epcot** è uno dei quattro parchi a tema che si trovano presso il *Walt Disney World Resort* a Bay Lake, Florida. Epcot rappresenta un'utopistica città del futuro. Il simbolo di questo parco è un'astronave chiamata *Spaceship Earth*[24]; essa ha la forma di una sfera geodetica di tipo (8, 8) (vedere la Figura 42). Le linee che sono evidenziate in questa figura rappresentano la proiezione sulla sfera geodetica dei segmenti che, nella tassellazione del piano, formano una linea che collega due vertici adiacenti dei triangoli grandi.

Essendo di tipo (8,8), si tratta di una sfera geodetica con $T$ = 192; dunque è formata da 192 · 20 = 3840 triangoli isosceli e ognuno di questi è la base di una piramide a base triangolare.

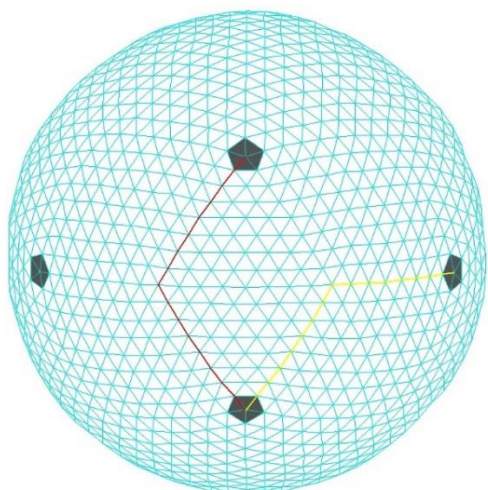

*Figura 42. Sfera geodetica di tipo (8, 8)*

---

[24] Buckminster Fuller chiamava la Terra "navicella spaziale".

Di conseguenza la cupola è composta da $3 \cdot 3840 = 11.520$ faccette. In realtà esse sono 11.324, poiché alcune di queste sono state eliminate dalla presenza dei supporti e delle porte.

## 8 – *Considerazioni finali*

La geometria delle cupole geodetiche ideate da Buckminster Fuller nasce dalle sue osservazioni sui fenomeni naturali: in Natura, i carichi in una struttura geometrica vengono distribuiti nel percorso più breve possibile (le linee geodetiche appunto) e i campi energetici "triangolano" lo spazio. Nella cupola geodetica, i triangoli (non uguali ma non troppo "diversi" fra loro) formano una griglia (energeticamente) compatta che garantisce la robustezza locale, mentre la disposizione dei loro lati lungo le geodetiche permette di distribuire gli sforzi locali sull'intera struttura, ottenendo, quindi, la massima efficienza col minor sforzo, sia a livello "energetico" che "materiale".

Il calcolo degli elementi geometrici che compongono una cupola geodetica (lunghezze delle aste, angoli delle facce, angoli diedri, etc.) una volta era lungo e difficile; sono state stilate varie tabelle da più autori, tenute sempre gelosamente custodite! Con l'avvento dei computer questi calcoli sono diventati più semplici ed oggi esistono diversi software realizzati proprio per determinare le misure delle varie entità geometriche di una cupola geodetica.

Anche se queste strutture, dopo un iniziale boom costruttivo, non hanno avuto quella diffusione auspicata da Buckminster Fuller, sono comunque oggetto di studio anche oggi e non solo in campo architettonico-ingegneristico ma anche in vari campi delle Scienze, dalla Biologia alla Chimica.

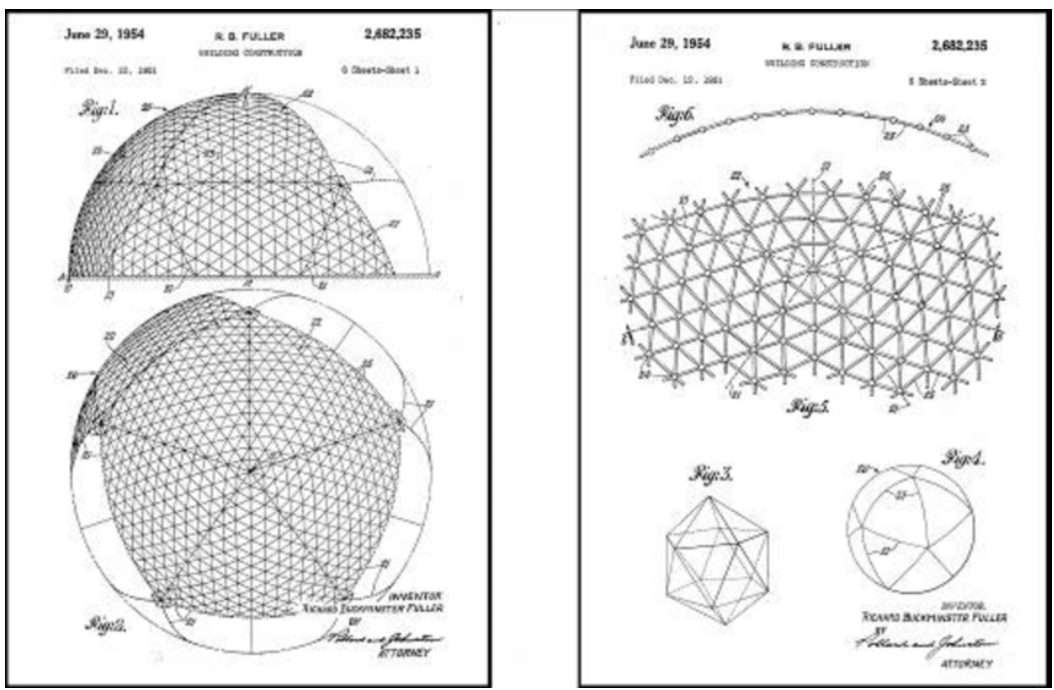

*Figura 43. Fogli 1 e 2 del brevetto di R. Buckminster Fuller US 2682235, 1954*

Se le prime sperimentazioni di Buckminster Fuller si sono basate su modelli in cui vengono riprodotti tutti i 31 cerchi massimi ricavati dai piani equatoriali delle simmetrie rotazionali di un icosaedro [Ke], ottenendo una struttura robusta ma esteticamente sgradevole perché molto disomogenea nelle lunghezze degli elementi e nel grado dei numerosi vertici, si è poi passati alla riduzione dei cerchi considerati ed alla successiva rimozione di alcune parti di questi togliendo spigoli dell'icosaedro di partenza. Le strutture così ottenute, con vertici di grado 5 o 6 e con elementi di lunghezze meno diverse fra loro, vengono poi eventualmente rinforzate introducendo elementi paralleli a quelli esistenti. Gli archi geodetici, che danno il nome alla cupola, sembrano "spariti". Nella descrizione dei 31 cerchi massimi dell'icosaedro, Hugh Kenner fa riferimento a tre tipologie (6+15+10) che nascono da rotazioni intorno a tre tipi di assi (una coppia di vertici opposti, i punti mediani di due spigoli opposti, i baricentri di due facce opposte). Ogni cerchio massimo può essere ricostruito a partire dal triangolo di Schwarz di una faccia dell'icosaedro proiettato sulla sfera circoscritta. Il triangolo di Schwarz di un triangolo equilatero si ottiene tracciando le tre altezze e considerando uno dei 6 triangoli retti (con angoli 30° e 60°) in cui risulta suddiviso il triangolo di partenza. Notiamo che, accoppiando un triangolo di Schwarz col suo simmetrico rispetto al cateto minore, si ottiene un triangolo isoscele che tassella il triangolo di partenza applicando solo rotazioni. Proiettando il triangolo di Schwarz sulla sfera, si ottiene un triangolo sferico che ne rappresenta 1/120 e che può tassellarla con riflessioni e rotazioni[25].

Si ha quindi la possibilità di tassellare il triangolo di Schwarz, anziché l'intera faccia dell'icosaedro, per poi riportare la tassellazione sulla sfera proiettando i vari vertici dal centro della sfera. Per una descrizione completa si veda ad esempio [Ki]. Lo stesso Walter Bauersfeld, dagli appunti ritrovati [Ga1], usa questo metodo per la costruzione della cupola di Jena.

Recentemente [H-N-G], per aumentare la regolarità e la "sfericità" di una cupola geodetica, è stata studiata una "stepping-projection", ovvero una proiezione a passi successivi. Partendo da una cupola con frequenza 2v, ogni faccia dell'icosaedro originale risulta suddivisa in 4 triangoli. Dividendo ognuno di questi triangoli con frequenza 2v e proiettando i nuovi vertici sulla sfera, si ottiene una cupola tipo 4v con elementi triangolari di forme meno "diverse" fra loro rispetto alla usuale costruzione 4v. Nell'articolo citato sono mostrati calcoli e tabelle relative a queste costruzioni.

## Ringraziamenti

## Bibliografia

---

[25] Proiettando sulla sfera il triangolo di Schwarz unito al suo simmetrico si può tassellare la sfera con 60 unità senza usare riflessioni. Ripetendo la costruzione del triangolo di Schwarz a partire da un tetraedro si ottiene un tassello che copre 1/24 della sfera mentre, partendo da un ottaedro, il triangolo sferico ottenuto ne copre 1/48.